\documentclass[11pt]{article}
\usepackage{amssymb}
\usepackage{mathrsfs}
\usepackage{latexsym,amsmath,amssymb,amsfonts,epsfig,graphicx,cite,psfrag}
\usepackage{eepic,color,colordvi,amscd}
\usepackage{comment}

\newtheorem{theo}{Theorem}[section]

\newtheorem{lem}[theo]{Lemma}

\def\qed{\hfill \rule{4pt}{7pt}}

\def\pf{\noindent {\it Proof. }}

\begin{document}

\title{The Kelmans-Seymour conjecture  II: 2-vertices in $K_4^-$}

\author{Dawei He\footnote{dhe9@math.gatech.edu; Partially supported by NSF grant through X. Yu},
Yan Wang\footnote{yanwang@gatech.edu; Partially supported by NSF grant through X. Yu}, 
Xingxing Yu\footnote{yu@math.gatech.edu; Partially supported by NSF grants DMS--1265564 and CNS--1443894} \\
\medskip\\
School of Mathematics\\
Georgia Institute of Technology\\
Atlanta, GA 30332-0160, USA}

\date{}

\maketitle

\begin{abstract}
We use $K_4^-$ to denote the graph obtained from $K_4$ by removing an edge, and use $TK_5$ to denote a subdivision 
of $K_5$.  Let $G$ be a 5-connected nonplanar graph and $\{x_1,x_2,y_1,y_2\}\subseteq V(G)$ such that 
$G[\{x_1,x_2,$ $y_1,y_2\}]\cong K_4^-$ with $y_1y_2\notin E(G)$. Let $w_1,w_2,w_3\in N(y_2)-\{x_1,x_2\}$ be 
distinct. We show that $G$ contains a $TK_5$ in which $y_2$ is not a branch vertex, or $G-y_2$ contains $K_4^-$, or $G$ has a special 
5-separation, or $G-\{y_2v:v\notin \{w_1,w_2,w_3,x_1,x_2\}\}$ contains $TK_5$. 
 

\end{abstract}

\newpage 

\section{Introduction}

We use notation and terminology from \cite{HWY15I}. In particular, for a graph $K$, we use $TK$ to denote  a \textit{subdivision} of $K$. 
The vertices in a $TK$ corresponding to the vertices of $K$ are its \textit{branch vertices}. 
Kelmans~\cite{Ke79} and, independently, Seymour~\cite{Se77} conjectured that every 5-connected nonplanar graph contains $TK_5$. 
In \cite{MY10,MY13}, this conjecture is shown to be true for  graphs containing $K_4^-$.

In \cite{HWY15I} we outline a strategy to prove the Kelmans-Seymour conjecture for graphs containing no $K_4^-$. 
Let $G$ be a 5-connected nonplanar graph containing no $K_4^-$. Then by a result of Kawarabayashi~\cite{Ka02}, $G$ contains an edge $e$ such that $G/e$ is 5-connected. If $G/e$ is planar, 
we can apply a discharging argument. So assume $G/e$ is not planar. Let $M$ be a maximal connected subgraph of $G$ such that $G/M$ is 5-connected and nonplanar. 
Let $z$ denote the vertex representing the contraction of $M$, and let $H=G/M$. Then one of the following holds:
\begin{itemize}
\item [(a)] $H$ contains a $K_4^-$ in which  $z$ is of degree 2.
\item [(b)] $H$ contains a $K_4^-$ in which  $z$ is of degree 3.
\item [(c)] $H$ does not contain $K_4^-$, and there exists $T\subseteq H$ such that $z\in V(T)$, $T\cong K_2$ or $T\cong K_3$, and $H/T$ is 5-connected and planar.
\item [(d)] $H$ does not contain  $K_4^-$, and for any $T\subseteq H$ with $z\in V(T)$ and $T\cong K_2$ or $T\cong K_3$, $H/T$ is not 5-connected.
\end{itemize}

In this paper, we deal with (a) by taking advantage of the $K_4^-$ containing $z$. We prove the following result, in which the vertex $y_2$ plays the role of $z$ above.

\begin{theo}\label{y_2}
Let $G$ be $\textsl{a}$ $5$-connected nonplanar graph and $\{x_1 , x_2 , y_1 , y_2\}\subseteq V(G)$ such that $G[\{x_1 , x_2 , y_1 , y_2\}]\cong K_4^-$ with $y_1y_2\notin E(G)$. 
Then one of the following holds:
\begin{itemize}
\item[$(i)$] $G$ contains $\textsl{a}$ $TK_5$ in which $y_2$ is not $\textsl{a}$ branch vertex.
\item[$(ii)$] $G-y_2$ contains $K_4^-$.
\item[$(iii)$] $G$ has $\textsl{a}$ $5$-separation $(G_1,G_2)$ such that $V(G_1\cap G_2)=\{y_2, a_1,a_2,a_3,a_4\}$, and  $G_2$ is 
the graph obtained from the edge-disjoint union of the $8$-cycle $a_1b_1a_2b_2a_3b_3a_4b_4\allowbreak a_1$ 
and the $4$-cycle $b_1b_2b_3b_4b_1$ by adding $y_2$ and the edges $y_2b_i$ for $i\in [4]$.
\item[$(iv)$] For $w_1,w_2,w_3\in N(y_2)-\{x_1,x_2\}$, $G-\{y_2v: v\notin \{w_1,w_2,w_3,x_1,x_2\}\}$ contains  $TK_5$.
\end{itemize}
\end{theo}

Note that when Theorem~\ref{y_2} is applied  later, $G$ will be a graph obtained from a 5-connected nonplanar graph by contracting a connected subgraph, and $y_2$ represents that contraction. 
So we need a $TK_5$ in $G$ to satisfy $(i)$ or $(iv)$ to produce a $TK_5$ in the original graph.  Note that $(ii)$ will not occur if the original graph is $K_4^-$-free. Moreover, 
if $(iii)$ occurs then we may apply Proposition 1.3 in \cite{HWY15I} to produce a $TK_5$ in the original graph.

The arguments used in this paper to prove Theorem~\ref{y_2} is similar to
those used in \cite{MY10, MY13}. Namely, we will find a substructure in the graph and use it to find the desired $TK_5$. However, since the $TK_5$ we are looking for must 
use certain special edges at $y_2$, the arguments here are more complicated and make heavy use of the option $(ii)$. 

We organize this paper as follows. In Section 2, we collect a few known results that will be used in the proof of Theorem~\ref{y_2}. 
We will produce an intermediate structure in $G$ which consists of eight special paths $X,Y,Z,A,B,C,P,Q$, see Figure~\ref{structure} (where $X$ is the path in bold and $Y,Z$ are not 
shown). In Section 3, we find the path $X$ in $G$ between $x_1$ and $x_2$ whose deletion  results in a graph satisfying certain connectivity requirement. 
In Section 4, we find the paths $Y,Z,A,B,C,P,Q$ in $G$. In Section 5, we use this structure to find the desired $TK_5$ for Theorem~\ref{y_2}.

\section{Previous results}

Let $G$ be a graph and $A\subseteq 
 V(G)$, and let $k$ be a positive
integer.
Let $[k] = \{1,2,...,k\}$. Let $C$ be a  cycle in $G$ with a fixed orientation (so that we can speak of clockwise and anticlockwise directions).
For two vertices $x,y \in V(C)$, $xCy$ denotes the subpath of $C$ from $x$ to $y$ in clockwise order. 
(If $x=y$ then $xCy$ denotes the path consisting of the single vertex $x$.)  
Recall from \cite{HWY15I} that $G$ is {\it $(k,A)$-connected} if, for any cut $T$ of
$G$ with $|T|< k$, every component of $G-T$ contains a vertex from
$A$. We say that $(G,A)$ is 
{\it plane} if $G$ is drawn in the plane with no crossing edges such that the vertices in $A$ are incident with the unbounded 
face of $G$. Moreover, for vertices $a_1,\ldots, a_k\in V(G)$, we say $(G, a_1,\ldots,a_k)$ is {\it plane} if 
$G$ is drawn in a closed disc in the plane with no crossing edges such that  
 $a_1,\ldots, a_k$ occur on the boundary of the disc in this cyclic order. We say that $(G,A)$ is {\it planar} if $G$ has a plane representation such that 
$(G,A)$ is plane. Similarly, $(G, a_1,\ldots,a_k)$ is {\it planar} if  $G$ has a plane representation such that $(G, a_1,\ldots,a_k)$ is plane.

In this section, we list a few known results that we need. We begin with a technical notion. 
A {\it 3-planar graph} $(G,{\cal A})$ consists
of a graph $G$ and  a collection ${\cal A}=\{A_1,\ldots,
A_k\}$ of pairwise 
disjoint subsets of $V(G)$ $($possibly ${\cal A}=\emptyset)$ such that
\begin{itemize} 
\addtolength{\baselineskip}{-1ex}

\item for distinct $i, j\in [k]$, $N(A_i)\cap A_j=\emptyset$, 

\item for $i\in [k]$,  $|N(A_i)|\le 3$, and 

\item  if $p(G,{\cal A})$ denotes the 
graph obtained from $G$ by (for each $i\in [k]$) deleting $A_i$ and
adding new edges joining every pair of distinct vertices in  $N(A_i)$,
then $p(G,{\cal A})$ can be drawn in a closed disc  with no 
crossing edges.
\end{itemize}
\noindent If, in addition,  $b_1, \ldots, b_n$ are vertices in $G$ such that $b_i\notin A_j$
for all $i\in [n]$ and $j\in [k]$, $p(G,A)$ can be drawn in a closed disc in the plane  with
no crossing edges, and $b_1,\ldots,b_n$ occur on the boundary of the disc in this cyclic order,
then we say that $(G,{\cal A}, b_1,\ldots,b_n)$ is
{\it 3-planar}. If there is no need to specify ${\cal A}$, we will simply say 
that $(G,b_1,\ldots,b_n)$ is {\it 3-planar}. 

It is easy to see that if  $(G,{\cal A}, b_1,\ldots,b_n)$ is
{\it 3-planar} and  $G$ is $(4,\{b_1,\ldots,b_n\})$-connected then ${\cal A}=\emptyset$ and  
$(G, b_1,\ldots,b_n)$ is planar. 

We can now state the following result of Seymour 
\cite{Se80}; equivalent versions can be found in \cite{CR79, Th80, Sh80}.

\begin{lem}\label{2path} 
Let $G$ be a graph and $s_1,s_2,t_1,t_2$ be distinct vertices of $G$. Then
exactly one of the following holds:
\begin{itemize}
\item [$(i)$] $G$ contains disjoint paths from $s_1$ to $t_1$ and from
  $s_2$ to $t_2$. 
\item [$(ii)$] $(G,s_1,s_2,t_1,t_2)$ is 3-planar.
\end{itemize}
\end{lem}

We also state a generalization of Lemma~\ref{2path},  which is a consequence of Theorems 2.3 and 2.4 in \cite{RS90}.

\begin{lem}\label{society}
Let $G$ be a graph,  $v_1,\ldots, v_n\in V(G)$ be distinct, and $n\ge 4$. Then 
exactly one of the following holds:
\begin{itemize}
\item [$(i)$] There exist $1\le i<j<k<l\le n$ such that $G$ contains disjoint paths from $v_i,v_j$ to $v_k, v_l$, respectively. 
\item [$(ii)$] $(G,v_1,v_2,\ldots,v_n)$ is 3-planar.
\end{itemize}
\end{lem}

The next result is Theorem 1.1 in \cite{HWY15I}. 

\begin{lem}\label{apexside1}
Let $G$ be $\textsl{a}$ $5$-connected nonplanar graph and let $(G_1, G_2)$ be $\textsl{a}$ $5$-separation in $G$. 
Suppose $|V(G_i)|\geq 7$ for $i\in [2]$, $a\in V(G_1\cap G_2)$, and $(G_2-a,V(G_1\cap G_2)-\{a\})$ is planar. Then one of the following holds:
\begin{itemize}
\item [$(i)$]  $G$ contains $\textsl{a}$ $TK_5$ in which $a$ is not $\textsl{a}$ branch vertex.
\item [$(ii)$] $G-a$ contains $K_4^-$.
\item [$(iii)$] $G$ has $\textsl{a}$ $5$-separation $(G_1',G_2')$ such that $V(G_1'\cap G_2')=\{a, a_1,a_2,a_3,a_4\}$, $G_1\subseteq G_1'$, and $G_2'$ is 
the graph obtained from the edge-disjoint union of the $8$-cycle $a_1b_1a_2b_2a_3b_3a_4\allowbreak b_4a_1$ 
and the  $4$-cycle $b_1b_2b_3b_4b_1$ by adding $a$ and the edges $ab_i$ for $i\in [4]$.
\end{itemize}
\end{lem}

Another result we need is Theorem 1.2 from \cite{HWY15I}. 

\begin{lem}
\label{5cut_triangle}
Let $G$ be \textsl{a} $5$-connected graph and $(G_1,G_2)$ be \textsl{a} $5$-separation in $G$. Suppose that $|V(G_i)|\ge 7$ for $i\in [2]$ and 
$G[V(G_1\cap G_2)]$ contains \textsl{a} triangle $aa_1a_2a$. Then one of the following holds: 
\begin{itemize}
\item [$(i)$]  $G$ contains \textsl{a} $TK_5$ in which $a$ is not \textsl{a} branch vertex.
\item [$(ii)$] $G-a$ contains $K_4^-$.
\item [$(iii)$] $G$ has \textsl{a} $5$-separation $(G_1',G_2')$ such that $V(G_1'\cap G_2')=\{a, a_1,a_2,a_3,a_4\}$ and $G_2'$ is 
the graph obtained from the edge-disjoint union of the $8$-cycle $a_1b_1a_2b_2a_3b_3a_4b_4\allowbreak a_1$ 
and the $4$-cycle $b_1b_2b_3b_4b_1$ by adding $a$ and the edges $ab_i$ for $i\in [4]$.
\item [$(iv)$] For any distinct $u_1,u_2,u_3\in N(a)-\{a_1,a_2\}$, $G-\{av: v \not\in \{a_1,a_2,u_1,u_2,u_3\}\}$ contains $TK_5$. 
\end{itemize}
\end{lem}

We also need Proposition 4.2 from \cite{HWY15I}. 

\begin{lem}\label{apex}
 Let $G$ be $\textsl{a}$ $5$-connected nonplanar graph and $a \in V (G)$ such that $G-a$ is
planar. Then one of the following holds:
\begin{itemize}
\item [$(i)$]  $G$ contains \textsl {a} $TK_5$ in which $a$ is not $\textsl{a}$ branch vertex.
\item [$(ii)$] $G- a$ contains $K_4^-$.
\item [$(iii)$] $G$ has $\textsl{a}$ $5$-separation $(G_1, G_2)$ such that $V (G_1 \cap G_2) = \{a, a_1, a_2, a_3, a_4\}$ and $G_2$ is the graph obtained 
from the edge-disjoint union of the 8-cycle $a_1b_1a_2b_2a_3b_3a_4b_4a_1$ and the 4-cycle $b_1b_2b_3b_4b_1$ by adding $a$ and the edges $ab_i$ for $i \in [4]$.
\end{itemize}
\end{lem}

We will make use of the following result of Perfect \cite{Pe68} on independent paths. A collection of paths in a graph are said to be {\it independent}
if no internal vertex of a path in this collection belongs to another path in the collection.

\begin{lem}
\label{perfect}
Let $G$ be a graph, $u\in V(G)$, and $A\subseteq V(G-u)$. Suppose there exist $k$ independent paths from $u$ to distinct $a_1,\ldots, a_k\in A$, respectively, 
and otherwise disjoint from $A$. Then for any $n\ge k$, if there exist $n$
independent paths  $P_1,\ldots, P_n$ in $G$ from $u$ to $n$ distinct vertices in $A$ and otherwise disjoint from $A$ then  $P_1,\ldots, P_n$ may be chosen so that $a_i\in V(P_i)$ for $i\in [k]$. 
\end{lem}

We will also use a result of Watkins and Mesner \cite{WM67} on cycles through three vertices. 
\begin{lem}
\label{Watkins}
Let $G$ be a $2$-connected graph and let $y_1, y_2, y_3$ be three distinct vertices of $G$. 
Then there is no cycle in $G$ containing $\{y_1, y_2, y_3\}$ if, and only if, one of the following statements holds: 
\begin{itemize}
\item [$(i)$] There exists a 2-cut $S$ in $G$ and there exist pairwise disjoint subgraphs $D_{y_i}$ of $G - S$, $i=1,2,3$, such that 
$y_i\in V(D_{y_i})$ and each $D_{y_i}$ is a union of components of $G - S$. 
\item [$(ii)$] There exist 2-cuts $S_{y_i}$ of $G$, $i=1,2,3$,  $z\in S_{y_1} \cap S_{y_2} \cap S_{y_3}$, and pairwise disjoint subgraphs $D_{y_i}$ of $G$, 
such that $y_i \in V(D_{y_i})$, each $D_{y_i}$ is a union of components of $G-S_{y_i}$, and $S_{y_1} - \{z\}, S_{y_2} - \{z\}, S_v - \{z\}$ are pairwise disjoint. 
\item [$(iii)$] There exist pairwise disjoint $2$-cuts $S_{y_i}$ in $G$, $i=1,2,3$, and pairwise disjoint subgraphs $D_{y_i}$ of $R - S_{y_i}$ such that $y_i \in V(D_{y_i})$, 
each $D_{y_i}$ is a union of components of $G - S_{y_i}$, and $G - V(D_{y_1} \cup D_{y_2} \cup D_{y_3})$ has precisely two components, 
each containing exactly one vertex from $S_{y_i}$ for $i\in [3]$.
\end{itemize}
\end{lem}

\section{Nonseparating paths}

Our first step for proving Theorem~\ref{y_2} is to find the path $X$ in $G$ (see Figure~\ref{structure}) whose removal does not affect connectivity too much. 

We need the concept of chain of blocks. Let $G$ be a graph and $\{u, v\}\subseteq V(G)$. We say that a sequence of blocks 
$B_1, \ldots, B_k$ in $G$ is a {\it chain of blocks} from $u$ to $v$ if  either $k=1$ and $u,v\in V(B_1)$ are distinct, or $k\ge 2$, 
$u\in V(B_1)-V(B_2)$, $v\in V(B_k)-V(B_{k-1})$, $|V(B_i)\cap V(B_{i+1})|=1$ for $i\in [k-1]$, 
and $V(B_i)\cap V(B_j)=\emptyset$ for any $i,j\in [k]$ with $|i-j|\ge 2$. For convenience, we also view this chain of blocks as $\bigcup_{i=1}^kB_i$, 
a subgraph of $G$.

The following result was implicit in \cite{CY03, KLY05}. Since it has not been stated and proved explicitly before, we 
include a proof.  We need the concept of a bridge. Let $G$ be a graph and $H$ a subgraph of $G$. Then an {\it $H$-bridge} of $G$
is a subgraph of $G$ that is either  induced by an edge of $G-E(H)$ with both ends in $V(H)$, or induced by the edges in some component of $G-H$ as well 
as those edges of $G$ from that component to $H$.

\begin{lem}\label{4A} 
Let $G$ be a graph and let $x_1, x_2, y_1, y_2\in V(G)$ be distinct such that $G$ is $(4, \{x_1, x_2, y_1, y_2\})$-connected. Suppose 
there exists a path $X$ in $G-x_1x_2$ from $x_1$ to $x_2$ such that $G-X$ contains a chain of blocks $B$ from $y_1$ to $y_2$. 
Then one of the following holds:
\begin{itemize}
 \item[$(i)$] There is a $4$-separation $(G_1, G_2)$ in $G$ such that $B+\{x_1,x_2\}\subseteq G_1$, $|V(G_2)|\geq 6$, and $(G_2, V(G_1\cap G_2))$ is planar.
 \item[$(ii)$] There exists an induced path $X'$ in $G-x_1x_2$ from $x_1$ to $x_2$ such that $G-X'$ is a chain of blocks from $y_1$ to $y_2$ and contains $B$.
\end{itemize}
\end{lem}

\pf Without loss of generality, we may assume that $X$ is induced in $G-x_1x_2$. 
We choose such $X$ that 
\begin{itemize}
\item [$ (1)$] $B$ is maximal,
\item [$ (2)$] the smallest size of a component of $G-X$ disjoint from $B$ (if exists)  is minimal, and 
\item [$(3)$]   the number of components of $G-X$ is minimal. 
\end{itemize}

We claim that $G-X$ is connected. For, suppose $G-X$ is not connected and let $D$ be a component of $G-X$ other than $B$ such that $|V(D)|$ 
is minimal. Let $u,v\in N(D)\cap V(X)$ such that $uXv$ is maximal. Since $G$ is $(4,\{x_1, x_2, y_1, y_2\})$-connected, $uXv-\{u,v\}$ contains 
a neighbor of some component of $G-X$ other than $D$. Let $Q$ be an induced path in $G[D+\{u,v\}]$ from $u$ to $v$, and let $X'$ be obtained from $X$ 
by replacing $uXv$ with $Q$. Then $B$ is contained in $B'$, the chain of blocks in $G-X'$ from $y_1$ to $y_2$. 
Moreover, either  the smallest size of a component of $G-X'$ disjoint from $B'$ is smaller than the smallest size of a component of 
$G-X$ disjoint from $B$, or the number of components of $G-X'$ is smaller than 
the number of components of $G-X$. This gives a contradiction to (1) or (2) or (3). Hence, $G-X$ is connected. 
 
If $G-X=B$, we are done with $X':=X$. So assume $G-X\neq B$. By (1), each $B$-bridge of $G-X$ has exactly one vertex in $B$. 
Thus, for each $B$-bridge $D$ of $G-X$, let $b_D\in V(D)\cap V(B)$ and $u_D,v_D\in N(D-b_D)\cap V(X)$ such that $u_DXv_D$ is maximal.

We now define a new graph ${\cal B}$ such that $V({\cal B})$ is the set of all $B$-bridges of $G-X$, and 
two $B$-bridges in $G-X$, $C$ and $D$, are adjacent if $u_CXv_C-\{u_C,v_C\}$ contains a neighbor of $D-b_D$ or  $u_DXv_D-\{u_D,v_D\}$ contains a neighbor of $C-b_C$. 
Let ${\cal D}$ be a component of ${\cal B}$. Then $\bigcup_{D\in V({\cal D})}u_DXv_D$ is a subpath of $X$. Let $S_{\cal D}$ be the union of $\{b_D:D\in V({\cal D})\}$ and 
the set of neighbors in $B$ of the internal vertices of $\bigcup_{D\in V({\cal D})}u_DXv_D$.

Suppose ${\cal B}$ has a component   ${\cal D}$ such that $|S_{\cal D}|\le 2$. 
Let $u,v\in V(X)$ such that $uXv=\bigcup_{D\in V({\cal D})}u_DXv_D$. Then 
$\{u,v\} \cup S_{\cal D}$ is a cut in $G$.
Since $G$ is $(4,\{x_1,x_2,y_1,y_2\})$-connected, $|S_{\cal D}| = 2$. 
So there is a $4$-separation $(G_1, G_2)$ in $G$ such that  $V(G_1\cap G_2) = \{u,v\}\cup S_{\cal D}$, $B+\{x_1,x_2\}\subseteq G_1$, and $D\subseteq G_2$ for $D\in V({\cal D})$.
Hence $|V(G_2)|\geq 6$.
If $G_2$ has disjoint paths $S_1,S_2$, with $S_1$ from  $u$ to $v$ and $S_2$ between the vertices in $S_{\cal D}$, then choose $S_1$ to be induced and 
let $X' = x_1Xu \cup S_1 \cup vXx_2$; now 
$B\cup S_2$ is contained in the chain of blocks in $G-X'$ from $y_1$ to $y_2$, contradicting (1). So no such two paths exist. 
Hence, by Lemma \ref{2path}, $(G_2, V(G_1\cap G_2))$ is planar and thus $(i)$ holds.

Therefore, we may assume that $|S_{\cal D}|\ge 3$ for any component ${\cal D}$ of ${\cal B}$. Hence, there exist a component ${\cal D}$ of ${\cal B}$ and
$D\in V({\cal D})$ with the following property:
$u_DXv_D-\{u_D,v_D\}$ contains vertices $w_1,w_2$ and $S_{\cal D}$ contains distinct vertices $b_1,b_2$ such that for each $i\in [2]$, $\{b_i,w_i\}$ is contained in 
a $(B\cup X)$-bridge of $G$ disjoint from $D-b_D$.    
Let $P$ denote an induced path in $G[D+\{u_{D},v_{D}\}]$ 
between $u_{D}$ and $v_{D}$, and let $X'$ be obtained from $X$ by
replacing $u_{D}Xv_{D}$ with $P$. Clearly, the chain of blocks in $G-X'$ from $y_1$ to $y_2$  
contains $B$ as well as a path from $b_1$ to $b_2$ and internally disjoint from $D\cup B$. This is a contradiction to (1).
\qed

\medskip

We now show that the conclusion of Theorem~\ref{y_2} holds or we can find a path $X$ in $G$ such that $y_1,y_2\notin V(X)$ and  $(G-y_2)-X$ is 2-connected.

\begin{lem}\label{reduction} 
Let $G$ be a $5$-connected nonplanar graph and let $x_1 , x_2 , y_1 , y_2\in V(G)$ be distinct such that $ G[\{x_1 , x_2 , y_1 , y_2 \}]\cong K_4^-$ 
with $y_1y_2\notin E(G)$. Then one of the following holds:
\begin{itemize}
 \item [$(i)$] $G$ contains \textsl{a} $TK_5$ in which $y_2$ is not \textsl{a} branch vertex.
\item [$(ii$)] $G-y_2$ contains $K_4^-$.
\item [$(iii)$] $G$ has \textsl{a} $5$-separation $(G_1,G_2)$ such that $V(G_1\cap G_2)=\{y_2, a_1,a_2,a_3,a_4\}$ and $G_2$ is 
the graph obtained from the edge-disjoint union of the $8$-cycle $a_1b_1a_2b_2a_3b_3a_4b_4\allowbreak a_1$ 
and the $4$-cycle $b_1b_2b_3b_4b_1$ by adding $y_2$ and the edges $y_2b_i$ for $i\in [4]$.
\item [$(iv)$]  For $w_1,w_2,w_3\in N(y_2)-\{x_1,x_2\}$,  $G-\{y_2v:v\notin \{w_1,w_2,w_3,x_1,x_2\}\}$ contains $TK_5$, or  
$G-x_1x_2$ has an induced path $X$ from $x_1$ to $x_2$ such that $y_1,y_2\notin V(X)$,  $w_1,w_2,w_3\in V(X)$, and $(G-y_2)-X$ is 2-connected.  
\end{itemize}
\end{lem} 
\pf First, we may assume that 
\begin{itemize}
\item [(1)]  $G-x_1x_2$ has an induced path $X$  from $x_1$ to $x_2$ such that $y_1,y_2\notin V(X)$ and $(G-y_2)-X$ is 2-connected.
\end{itemize} 
To see this, let $z\in N(y_1)-\{x_1,x_2\}$. Since $G$ is 5-connected, $(G-x_1x_2)-\{y_1,y_2,z\}$ has a path $X$ from $x_1$ to $x_2$. Thus,  we may apply 
Lemma~\ref{4A} to $G-y_2$, $X$ and $B=y_1z$.

Suppose $(i)$ of Lemma~\ref{4A} holds. Then 
$G$ has a $5$-separation $(G_1, G_2)$ such that $y_2\in V(G_1\cap G_2)$, $\{x_1,x_2,y_1,z\}\subseteq V(G_1)$ and $y_1z\in E(G_1)$, $|V(G_2)|\geq 7$, and $(G_2-y_2, V(G_1\cap G_2)-\{y_2\})$ is planar.  
If $|V(G_1)|\ge 7$ then, by Lemma~\ref{apexside1}, $(i)$ or $(ii)$ or $(iii)$  holds. If $|V(G_1)|=5$ then $G_1-y_2$ has a $K_4^-$ or $G-y_2$ is planar; hence, $(ii)$ holds in the former case, 
 and  $(i)$ or $(ii)$ or $(iii)$  holds in the latter case by Lemma~\ref{apex}. 
Thus we may assume that $|V(G_1)|=6$.  Let $v\in V(G_1-G_2)$. Then $v\ne y_2$. Since $G$ is 5-connected, $v$ must be adjacent to all vertices in 
$V(G_1\cap G_2)$. Thus, $v\ne y_1$ as $y_1y_2\notin E(G)$. Now $|V(G_1\cap G_2)\cap \{x_1,x_2,z\}|\ge 2$. Therefore, $G[\{v,y_1\}\cup (V(G_1\cap G_2)\cap \{x_1,x_2,z\})]$ 
contains $K_4^-$; so $(ii)$ holds. 

 So we may assume that $(ii)$ of Lemma~\ref{4A} holds. Then 
$(G-y_2)-x_1x_2$ has an induced path, also denoted by $X$,  from $x_1$ to $x_2$ such that $(G-y_2)-X$ is a chain of  blocks from $y_1$ to $z$. Since 
$zy_1\in E(G)$, $(G-y_2)-X$ is in fact a block. 
If $V((G-y_2)-X)=\{y_1,z\}$ then, since $G$ is  $5$-connected and $X$ is induced in $(G-y_2)-x_1x_2$, 
$G[\{x_1, x_2, z, y_1\}]\cong K_4$; so $(ii)$ holds. This completes the proof of (1). 

\medskip

We wish to prove $(iv)$. So let $w_1,w_2,w_3\in N(y_2)-\{x_1,x_2\}$ and assume that  
$$G':=G-\{y_2v:v\notin \{w_1,w_2,w_3,x_1,x_2\}\}$$ 
does not contain $TK_5$.  We may assume that

\begin{itemize}
\item [(2)] $w_1,w_2,w_3\notin V(X)$. 
\end{itemize}
For, suppose not. If $w_1,w_2,w_3\in V(X)$ then $(iv)$ holds. So, without loss of generality, we may assume 
$w_1\in V(X)-\{x_1, x_2\}$ and $w_2\in V(G-X)$. Since $X$ is induced in $G-x_1x_2$ and $G$ is $5$-connected, 
$(G-y_2)-(X-w_1)$ is 2-connected and, hence, contains independent paths 
$P_1, P_2$ from $y_1$ to $w_1, w_2$, respectively. Then $w_1Xx_1\cup w_1Xx_2\cup w_1y_2\cup P_1\cup (y_2w_2\cup P_2)\cup G[\{x_1, x_2, y_1, y_2\}]$ 
is a $TK_5$ in $G'$ with branch vertices $w_1,x_1,x_2,y_1,y_2$, a contradiction.

\begin{itemize}
\item [(3)]
For any $u\in V(x_1Xx_2)-\{x_1, x_2\}$, $\{u,y_1,y_2\}$ is not contained in any cycle in $G'-(X-u)$.
\end{itemize}
For, suppose there exists  $u\in V(x_1Xx_2)-\{x_1, x_2\}$ such that $\{u,y_1,y_2\}$ is contained in  a cycle $C$ in $G'-(X-u)$.
Then $uXx_1\cup uXx_2\cup C\cup G[\{x_1, x_2, y_1, y_2\}]$ is a $TK_5$ in $G'$ with branch vertices $u,x_1,x_2,y_1,y_2$, a contradiction. So we have (3). 

\medskip

Let $y_3\in V(X)$ such that $y_3x_2\in E(X)$, and let $H:= G'-(X-y_3)$. Note that $H$ is 2-connected. By (3), no cycle in $H$ contains 
$\{y_1, y_2,y_3\}$. Thus,  we  apply Lemma~\ref{Watkins} to $H$. In order to treat simultaneously the  three cases in the conclusion of Lemma~\ref{Watkins}, 
we introduce some notation. Let $S_{y_i}=\{a_i,b_i\}$ for $i\in [3]$, such that 
if  Lemma~\ref{Watkins}$(i)$ occurs we let  $a_1=a_2=a_3$, $b_1=b_2=b_3$, and $S_{y_i}=S$ for $i\in [3]$; if Lemma~\ref{Watkins}$(ii)$ occurs
then $a_1=a_2=a_3$; and if  Lemma~\ref{Watkins}$(iii)$ then $\{a_1,a_2,a_3\}$ and 
$\{b_1,b_2,b_3\}$ belong to different components of $H-V(D_{y_1}\cup D_{y_2}\cup D_{y_3})$. 
If  Lemma~\ref{Watkins}$(ii)$ or  Lemma~\ref{Watkins}$(iii)$ occurs then let $B_a, B_b$ denote the components of 
$H-V(D_{y_1}\cup D_{y_2}\cup D_{y_3})$ such that for $i\in [3]$ 
$a_i\in V(B_a)$ and $b_i\in V(B_b)$. Note that $B_a=B_b$ is possible, but only if Lemma~\ref{Watkins}$(ii)$ occurs.

For convenience, let $D_i':=G'[D_{y_i}+\{a_i,b_i\}]$ for $i\in [3]$. We choose the cuts $S_{y_i}$ so that 
\begin{itemize}
\item [(4)] $D_1'\cup D_2'\cup D_3'$ is maximal. 
\end{itemize}
Since $H$ is 2-connected,  $D_i'$, for each $i\in [3]$, contains a path $Y_i$ from $a_i$ to $b_i$ and through $y_i$. 
In addition, since $(G-y_2)-X$ is 2-connected,  for any $v\in V(D_3')-\{a_3,b_3,y_3\}$, $D_3'-y_3$ contains a path from $a_3$ to $b_3$ through $v$.

\begin{itemize}
\item [(5)]
If $B_a\cap B_b=\emptyset$ then $|V(B_a)|=1$ or $B_a$ is 2-connected, and $|V(B_b)|=1$ or $B_b$ is 2-connected. 
If $B_a\cap B_b\ne \emptyset$ then $B_a=B_b$ and $B_a-a_3$ is 2-connected.
\end{itemize}
First, suppose $B_a\cap B_b=\emptyset$. By symmetry, we only prove the claim for $B_a$. Suppose  $|V(B_a)|>1$ and  
$B_a$ is not 2-connected. Then $B_a$ has a separation $(B_1,B_2)$ such that $|V(B_1\cap B_2)|\le 1$. Since $H$ is 2-connected, 
$|V(B_1\cap B_2)|=1$   and, for some permutation $ijk$ of $[3]$,
$a_i\in V(B_1)-V(B_2)$ and $a_j,a_k\in V(B_2)$. 
Replacing $S_{y_i},D_i'$ by $V(B_1\cap B_2)\cup \{b_i\},D_i'\cup B_1$, respectively, while  keeping $S_{y_j},D_j',S_{y_k},D_k'$ unchanged,  we derive a contradiction to (4). 

Now assume   $B_a\cap B_b\ne \emptyset$. Then $B_a=B_b$ by definition, and $a_1=a_2=a_3$ by our assumption above. 
Suppose  $B_a - a_3$ is not 2-connected. Then $B_a$ has a 2-separation $(B_1,B_2)$ with $a_3\in V(B_1\cap B_2)$. First, suppose for some 
permutation $ijk$ of $[3]$, $b_i\in V(B_1)-V(B_2)$ and $b_j,b_k\in V(B_2)$. 
Then replacing $S_{y_i},D_i'$ by $V(B_1\cap B_2),D_i'\cup B_1$, respectively, while keeping $S_{y_j},D_j',S_{y_k},D_k'$ unchanged,  we derive a contradiction to (4).
Therefore, we may assume  $\{b_1, b_2, b_3\}\subseteq V(B_1)$. Since $G$ is $5$-connected, there exists $rr'\in E(G)$ such that $r\in V(X)-\{y_3,x_2\}$
and $r'\in V(B_2-B_1)$. Let $R$ be  a path $B_2-(B_1-a_3)$ from $a_3$ to $r'$, and $R'$ a path in $B_1-B_2$ from $b_1$ to $b_2$. 
Then $(R\cup r'r\cup rXx_1)\cup (a_3Y_3y_3\cup y_3x_2)\cup a_3Y_1y_1 \cup a_3Y_2y_2\cup 
(y_1Y_1b_1\cup R'\cup b_2Y_2y_2)\cup G[\{x_1, x_2, y_1, y_2\}]$ is a $TK_5$ in $G'$ with branch vertices $a_3,x_1,x_2,y_1,y_2$, a contradiction.

\begin{itemize}
\item [(6)]  $D_{y_i}$ is connected for $i\in [3]$.
\end{itemize}
Suppose $D_{y_i}$ is not connected for some $i\in [3]$, and let $D$ be a component of $D_{y_i}$ not containing $y_i$. Since $G$ is $5$-connected, 
there exists $rr'\in E(G)$ such that $r\in V(X)-\{x_2,y_3\}$ and $r'\in V(D)$. 

Let $R$ be a path in $G[D+a_i]$ from $a_i$ to $r'$, and $R'$ a path from $b_1$ to $b_2$ in $B_b-a_3$. 
By (5), let $A_1,A_2,A_3$ be independent paths in $B_a$ from $a_i$ to $a_1,a_2,a_3$, respectively.   
Then $(R\cup r'r\cup rXx_1)\cup (A_1\cup a_1Y_1y_1)\cup (A_2\cup a_2Y_2y_2)\cup (A_3\cup a_3Y_3y_3\cup y_3x_2)\cup 
(y_1Y_1b_1\cup R'\cup b_2Y_2y_2)\cup G[\{x_1, x_2, y_1, y_2\}]$ is 
a $TK_5$ in $G'$ with branch vertices $a_i,x_1,x_2,y_1,y_2$, a contradiction.

\begin{itemize}
\item [(7)] If $a_1=a_2=a_3$ then $N(a_3)\cap V(X-\{x_2,y_3\})=\emptyset$. 
\end{itemize}
For, suppose $a_1=a_2=a_3$ and there exists $u\in N(a_3)\cap V(X-\{x_2,y_3\})$. Let $Q$ be a path in $B_b-a_3$ between $b_1$ and $b_2$, and let $P$ be a path in $D_3'-b_3$ from $a_3$ to $y_3$. 
Then $(a_3u\cup uXx_1)\cup (P\cup y_3x_2)\cup a_3Y_1y_1\cup a_3Y_2y_2\cup (y_1Y_1b_1\cup Q\cup b_2Y_2y_2)\cup G[\{x_1,x_2,y_1,y_2\}]$ 
is a $TK_5$ in $G'$ with branch vertices $a_3,x_1,x_2,y_1,y_2$, a contradiction.

\medskip
We may assume that
\begin{itemize}
\item [(8)] there exists  $u\in V(X)-\{x_1, x_2,y_3\}$ such that $N(u)-\{y_2\} \not\subseteq V(X\cup D_3')$. 
\end{itemize}
 For, suppose no such vertex exists. Then  $G$ has a 5-separation $(G_1,G_2)$ such that 
$V(G_1\cap G_2)=\{ a_3, b_3, x_1, x_2, y_2\}$, $X\cup D_3'\subseteq G_1$, and $D_1'\cup D_2'\cup B_a\cup B_b \subseteq G_2$. 
Clearly, $|V(G_2)|\ge 7$ since $|N(y_1)| \geq 5$ and $y_1y_2\notin E(G)$. If $|V(G_1)|\ge 7$ then, by Lemma~\ref{5cut_triangle}, 
$(i)$ or $(ii)$ or $(iii)$ or $(iv)$ holds. So we may assume $|V(G_1)|=6$. Then $X=x_1y_3x_2$ and $V(D_{y_3})=\{y_3\}$. 
Hence, $G[\{x_1,x_2,y_1,y_3\}]\cong K_4^-$; so  $(ii)$ holds.

\begin{itemize}
\item [(9)] For all $u\in V(X)-\{x_1,x_2,y_3\}$ with $N(u)-\{y_2\}\not\subseteq V(X\cup D_3')$, $N(u)\cap V(D_3'-y_3)=\emptyset$. 
\end{itemize}
For, suppose there exist $u\in V(X)-\{x_1,x_2,y_3\}$, $u_1\in (N(u)-\{y_2\})-V(X\cup D_3')$, and $u_2\in N(u)\cap V(D_3'-y_3)$. 
Recall (see before (5)) that there is a path $Y_3'$ in $D_3'-y_3$ from $a_3$ to $b_3$ through $u_2$. 

Suppose $u_1\in V(D_{y_i})$ for some $i\in [2]$. Then $D_i'-b_i$ (or $D_i'-a_i$) has a path $Y_i'$ from $u_1$ to $a_i$ (or $b_i$) through $y_i$. 
If $Y_i'$ ends at $a_i$ then let $P_a$, $P_b$ be disjoint paths in $B_a\cup B_b$ from $a_1$, $b_3$ to $a_2$, $b_{3-i}$, respectively; now 
$Y_i'\cup P_a\cup Y_{3-i}\cup P_b\cup b_3Y_3'u_2\cup u_2uu_1$ is a cycle in $G'-(X-u)$ containing  $\{u,y_1, y_2\}$, contradicting (3). 
So $Y_i'$ ends at $b_i$.  Let $P_b$, $P_a$ be disjoint paths in $B_a\cup B_b$ from $b_1$, $a_{3-i}$ to $b_2$, $a_3$, respectively. 
Then $Y_i'\cup P_b\cup Y_{3-i}\cup P_a\cup a_3Y_3'u_2\cup u_2uu_1$ is a cycle in $G'-(X-u)$ containing $\{u,y_1, y_2\}$,  contradicting  (3).

Thus, $u_1\in V(B_a\cup B_b)$. By symmetry and (7), assume $u_1\in
V(B_b)$. Note that $u_1\notin \{a_3,b_3\}$ (by the choice of $u_1$) and  $B_b-a_3$ is $2$-connected (by (5)). 
Hence, $B_b-a_3$ has disjoint paths $Q_1, Q_2$  from $\{u_1, b_3\}$ to $\{b_1, b_2\}$. By symmetry between $b_1$ and $b_2$, 
we may assume $Q_1$ is between $u_1$ and $b_1$ and $Q_2$ is between $b_3$ and $b_2$. Let $P$ be a path in $B_a$ from $a_1$ to $a_2$
(which is trivial if $|V(B_a)|=1$).  
Then $Q_1\cup u_1uu_2\cup u_2Y_3'b_3\cup Q_2\cup Y_2\cup P\cup Y_1$ is a cycle in $G'-(X-u)$ containing $\{y_1, y_2, u\}$, contradicting (3).

\begin{itemize}
\item [(10)] For any $u\in V(X)-\{x_1,x_2,y_3\}$ with $N(u)-\{y_2\}\not\subseteq V(X\cup D_3')$, there exists $i\in [2]$ such that
$N(u)-\{y_2\}\subseteq V(D'_i)$ and $\{a_i,b_i\}\not\subseteq N(u)$.
\end{itemize}
To see this, let $u_1,u_2\in  (N(u)-\{y_2\})-V(X\cup D_3')$ be distinct, which exist by (9) (and since $X$ is induced in $G'-x_1x_2$). 
Suppose we may choose such $u_1,u_2$ so that $\{u_1,u_2\}\not\subseteq V(D_i')$ for $i\in [2]$.

We claim that  $\{u_1,u_2\}\not\subseteq V(B_a)$ and   $\{u_1,u_2\}\not\subseteq  V(B_b)$.   
Recall that if $B_a\cap B_b\ne \emptyset$ then $B_a=B_b$ and if $B_a\cap B_b
=\emptyset$ then there is symmetry between $B_a$ and $B_b$. So if the claim fails we may assume that  $u_1,u_2\in V(B_b)$. 
Then by (5), $B_b-a_3$ is 2-connected; so $B_b-a_3$ contains disjoint paths $Q_1,Q_2$ from $\{u_1,u_2\}$ to $\{b_1,b_2\}$.
If $B_a=B_b$, let $P = a_3$.
If $B_a \cap B_b = \emptyset$, then let $P$ be a path in $B_a$ from $a_1$ to $a_2$.
Now $Q_1\cup u_1uu_2\cup Q_2 \cup Y_1 \cup P \cup Y_2 $ is a cycle in $G'-(X-u)$ 
containing $\{u,y_1,y_2\}$, contradicting (3).

Next, we show that $\{a_i,b_i\}\not\subseteq N(u)$ for $i\in [2]$.
For, suppose $u_1=a_i$ and $u_2=b_i$ for some $i \in [2]$. Then, since $\{u_1,u_2\}\cap \{a_3,b_3\}=\emptyset$, 
$|V(B_a)|\ge 2$ and $|V(B_b)|\ge 2$. By (5),  let $P_1,P_2$ be independent  paths in $B_a$ from $a_i$ to $a_{3-i},a_3$, respectively, 
and $Q_1,Q_2$ be independent paths in $B_b$ from $b_i$ to $b_{3-i},b_3$, respectively. 
Now $ua_i \cup ub_i \cup a_iY_iy_i \cup b_iY_iy_i \cup (y_ix_1\cup x_1Xu)\cup (P_1 \cup Y_{3-i} \cup Q_1) \cup (P_2 \cup a_3Y_3y_3) \cup (Q_2 \cup b_3Y_3y_3) 
\cup uXy_3 \cup y_ix_2y_3$ is a $TK_5$ in $G'$ with branch vertices $a_i,b_i,u,y_i,y_3$, a contradiction. 

Suppose $u_1\in V(B_a-a_3)$ and $u_2\in V(B_b-b_3)$. Then  $|V(B_a)|\ge 2$ and $|V(B_b)|\ge 2$. Let $Y_3'$ be a path in $D_3'-y_3$ from $a_3$ to $b_3$. 
First,  assume that $u_1\in \{a_1,a_2\}$ or $u_2\in \{b_1,b_2\}$. By symmetry, we may assume $u_1=a_1$. 
So $u_2\ne b_1$. By (5), $B_a-a_1$ contains a path $P$ from $a_2$ to $a_3$, and $B_b$ contains disjoint paths $Q_1,Q_2$ 
from $\{b_2,b_3\}$ to $b_1,u_2$, respectively. Then $Y_1\cup Q_1 \cup Y_2\cup P\cup Y_3'\cup Q_2\cup u_1uu_2$ is a cycle in $G'-(X-u)$ 
containing $\{u,y_1,y_2\}$,  contradicting (3). 
So $u_1\notin \{a_1,a_2\}$ and  $u_2\notin \{b_1,b_2\}$. Then by (5) and symmetry, we may assume that 
$B_a$ contains disjoint paths $P_1,P_2$ from $u_1,a_3$ to $a_1,a_2$, respectively. By (5) again, $B_b$ contains disjoint paths $Q_1,Q_2$ from 
$b_1,u_2$, respectively to $\{b_2,b_3\}$.
Now  $P_1 \cup Y_1 \cup Q_1 \cup Y_2 \cup P_2 \cup Y_3' \cup Q_2\cup u_2uu_1$ is a cycle in $G'-(X-u)$ 
containing $\{u,y_1,y_2\}$, contradicting (3). 

Therefore, we may assume $u_1\in V(D_{y_i})$ for some $i\in [2]$. 
By symmetry, we may assume that $u_1\in V(D_{y_1})$ and $D_1'-a_1$ contains a path $R_1$ from $u_1$ to $b_1$ and through $y_1$. 
Then $u_2\notin  V(D_1')$ as we assumed $\{u_1,u_2\}\not\subseteq V(D_1')$. 

Suppose $u_2\in V(D_{y_{2}})$. 
If $D_2'-a_2$ contains a path $R_2$ from $u_2$ to $b_2$ through $y_2$ then let $Q$ be a path in $B_b$ from $b_1$ to $b_2$; now $R_1\cup Q\cup R_2\cup u_2uu_1$ is a cycle in 
$G'-(X-u)$ containing $\{u,y_1,y_2\}$, contradicting (3). 
So $D_2'-b_2$ contains a path $R_2$ from $u_2$ to $a_2$ and through $y_2$. Now let $P$ be a path in $B_a$ from $a_2$ to $a_3$, $Q$ be a path in $B_b-a_3$ from $b_1$ to $b_3$. 
Let $Y_3'$ be a path in $D_3'-y_3$ from $a_3$ to $b_3$. Then 
$R_1\cup Q\cup Y_3'\cup P\cup R_2\cup u_2uu_1$  is a cycle in 
$G'-(X-u)$ containing $\{u,y_1,y_2\}$, contradicting (3). 

Finally, assume $u_2\in V(B_a\cup B_b)$. If $u_2\in V(B_b)$ then, by (5), let $Q_1,Q_2$ be disjoint paths in $B_b-a_3$ from $b_1,u_2$, respectively,  to 
$\{b_2, b_3\}$, and let $P$ be a path in $B_a$ from $a_2$ to $a_3$; now $u_2uu_1\cup R_1\cup Q_1\cup Q_2\cup Y_2\cup P\cup Y_3'$ is a cycle in 
$G'-(X-u)$ containing $\{u,y_1,y_2\}$, contradicting (3). So $u_2\notin V(B_b)$ and $u_2\in V(B_a-a_1)$; hence $B_a\cap B_b=\emptyset$.  Let 
$P$ be a path in $B_a$ from $u_2$ to $a_2$ and $Q$ be a path in $B_b$ from $b_1$ to $b_2$. Then $u_2uu_1\cup R_1\cup Q\cup Y_2\cup P$ is a cycle in 
$G'-(X-u)$ containing $\{u,y_1,y_2\}$, contradicting (3).  This completes the proof of (10). 
  
\medskip
By (10) and by symmetry, let $u\in V(X)-\{x_1,x_2,y_3\}$ and $u_1,u_2\in N(u)$ such that $u_1\in V(D_{y_1})$ and $u_2\in V(D_1')$. 
If $G[D_1'+u]$ contains independent paths $R_1,R_2$ from $u$ to $a_1,b_1$, respectively, such that $y_1\in V(R_1\cup R_2)$, 
then let $P$ be a path in $B_a$ between $a_1$ and $a_2$ and $Q$ be a path in $B_b-a_3$ between $b_1$ and $b_2$; now $R_1\cup P\cup Y_2\cup Q \cup R_2$ is a 
cycle in $G'-(X-u)$ containing $\{u,y_1,y_2\}$, contradicting (3). 
So such paths do not exist. Then  in the 2-connected graph $D_1^*:=G[D_1'+u]+\{c,ca_1,cb_1\}$ (by adding a new vertex $c$), 
there is no  cycle containing  $\{c,u,y_1\}$. Hence, by Lemma~\ref{Watkins}, $D_1^*$ has a 2-cut $T$ separating $y_1$ from $\{u,c\}$, and 
$T\cap \{u,c\}=\emptyset$. 

We choose $u,u_1,u_2$ and $T$ so that the $T$-bridge of $D_1^*$ containing $y_1$, denoted  $B$,  is minimal. Then $B-T$ contains no neighbor of $X-\{x_1,x_2\}$. 
Hence, $G$ has a 5-separation $(G_1,G_2)$ such that $V(G_1\cap G_2)= \{x_1,x_2,y_2\}\cup V(T)$, $B\subseteq G_1$, and $X\cup D_2'\cup D_3'\subseteq 
G_2$. Clearly, $|V(G_2)|\ge 7$. Since $y_1y_2\notin E(G)$ and $G$ is 5-connected,  
$|V(G_1)|\ge 7$. So  $(i)$ or $(ii)$ or $(iii)$ or $(iv)$  holds by Lemma~\ref{5cut_triangle}. \qed

\section{An intermediate substructure}

By Lemma~\ref{reduction}, to prove Theorem~\ref{y_2} it suffices to deal with the second part of $(iv)$ of Lemma~\ref{reduction}. 
Thus, let $G$ be a 5-connected nonplanar graph and $x_1,x_2,y_1,y_2\in V(G)$ be distinct such that $G[\{x_1,x_2,y_1,y_2\}]\cong K_4^-$ with $y_1y_2\notin E(G)$, 
let $w_1,w_2,w_3\in N(y_2)-\{x_1,x_2\}$ be distinct, and let $P$ be an induced path in 
 $G-x_1x_2$ from $x_1$ to $x_2$ such that $y_1,y_2\notin V(P)$, $w_1,w_2,w_3\in V(P)$, and  $(G-y_2)-P$ is 2-connected. 

Without loss of generality, assume $x_1,w_1,w_2,w_3,x_2$ occur on $P$ in order. Let $$X:=x_1Pw_1\cup w_1y_2w_3\cup w_3Px_2,$$ and let 
$$G':=G-\{y_2v:v\notin \{w_1,w_2,w_3,x_1,x_2\}\}.$$ Then $X$ is an induced path in $G'-x_1x_2$, $y_1\notin V(X)$, and $G'-X$ is 2-connected. 
For convenience, we record this situation by calling $(G, X, x_1, x_2, y_1, y_2, w_1,w_2,w_3)$ a {\it $9$-tuple}.

In this section, we obtain a substructure of $G'$ 
in terms of $X$ and seven additional paths $A,B,C, P,Q,Y,Z$ in $G'$. See Figure~\ref{structure}, 
where $X$ is the path in boldface and $Y,Z$ are not shown. 
First, we find  two special paths $Y,Z$ in $G'$ with Lemma~\ref{YZ} below. 
We will then use Lemma~\ref{ABC} to find the paths $A,B,C$, and 
use Lemma~\ref{PQ} to find the paths $P$ and $Q$. 
In the next section, we will use this substructure to find the desired $TK_5$ in $G$ or $G'$.

\begin{lem}\label{YZ}
Let $(G,X,x_1,x_2,y_1,y_2, w_1,w_2,w_3)$ be a $9$-tuple. Then one of the following holds:
\begin{itemize}
\item [$(i)$] $G$ contains $TK_5$ in which $y_2$ is not a branch vertex, or $G'$ contains $TK_5$.
\item [$(ii)$] $G-y_2$ contains $K_4^-$.
\item [$(iii)$] $G$ has $\textsl{a}$ $5$-separation $(G_1,G_2)$ such that $V(G_1\cap G_2)=\{y_2, a_1,a_2,a_3,a_4\}$, $G_2$ is 
the graph obtained from the edge-disjoint union of the $8$-cycle $a_1b_1a_2b_2a_3b_3a_4b_4a_1$ 
and the $4$-cycle $b_1b_2b_3b_4b_1$ by adding $y_2$ and the edges $y_2b_i$ for $i\in [4]$.
\item [(iv)] There exist $z_1\in V(x_1Xy_2)-\{x_1,y_2\}$, $z_2\in V(x_2Xy_2)-\{x_2, y_2\}$ 
such that $H:=G'-(V(X-\{y_2, z_1,z_2\})\cup E(X))$ has disjoint paths $Y,Z$ from $y_1,z_1$ to $y_2, z_2$, respectively.
\end{itemize}
\end{lem}

\pf Let $K$ be the graph obtained from $G-\{x_1, x_2,y_2\}$ by contracting $x_iXy_2-\{x_i, y_2\}$ to the new vertex $u_i$, for $i\in [2]$. 
Note that $K$ is 2-connected; since $G$ is 5-connected, $X$ is induced in $G'-x_1x_2$, and $G-X$ is 2-connected. We may assume that 

\begin{itemize}
\item [(1)] there exists a collection ${\cal A}$ of subsets of $V(K)-\{u_1,u_2,w_2,y_1\}$ such that $(K,{\cal A}, u_1,y_1,u_2,w_2)$ is 3-planar. 
\end{itemize}
For, suppose this is not the case. Then by Lemma~\ref{2path},  $K$ contains disjoint paths, say $Y,U$, from $y_1, u_1$ to $w_2, u_2$, respectively. 
Let $v_i$ denote the neighbor of $u_i$ in the path $U$, and let $z_i\in V(x_iXy_2)-\{x_i,y_2\}$ be a neighbor of $v_i$ in $G$. 
Then $Z:=(U-\{u_1, u_2\})+\{z_1, z_2, z_1v_1, z_2v_2\}$ is a path between $z_1$ and $z_2$. Now $Y+\{y_2, y_2w_2\}, Z$ are the desired paths for $(iv)$. 
So we may assume (1). 

\medskip
Since $G-X$ is 2-connected, $|N_K(A)\cap \{u_1,u_2, w_2\}|\le 1$ for all $A\in {\cal A}$. 
Let $p(K,{\cal A})$ be the graph obtained from $K$ by (for each $A\in \cal{A}$) deleting $A$ 
and adding new edges joining every pair of distinct vertices in $N_{K}(A)$. Since $G$ is 5-connected and $G-X$ is
2-connected, we may assume that $p(K,{\cal A})-\{u_1,u_2\}$ is a 2-connected plane
graph, and for each $A\in {\cal A}$ with $N_K(A)\cap \{u_1,u_2\}\ne \emptyset$ the edge joining vertices of $N_K(A)-\{u_1,u_2\}$
occur on the outer cycle $D$ of $p(K, {\cal A})-\{u_1, u_2\}$. Note that  $y_1,w_2\in V(D)$.

Let $t_1\in V(D)$ with $t_1Dy_1$ minimal such that $u_1t_1\in E(p(K,{\cal A}))$; and let $t_2\in V(D)$ with $y_1Dt_2$ 
minimal such that $u_2t_2\in E(p(K,{\cal A}))$. (So $t_1,y_1,t_2,w_2$ occur on $D$ in clockwise order.) 
Since $K$ is 2-connected and $X$ is induced in $G'-x_1x_2$, 
there exist $z_1\in V(x_1Xy_2)-\{x_1, y_2\}$ and independent paths $R_1, R_1'$ in $G$ 
from $z_1$ to $D$ and internally disjoint from $V(p(K,{\cal A}))\cup V(X)$, such that $R_1$ ends at $t_1$ and $R_1'$ ends at some vertex $t_1'\ne t_1$, and  
$w_2,t_1',t_1,y_1$ occur on $D$ in clockwise order. Similarly, there exist $z_2\in V(x_2Xy_2)-\{x_2, y_2\}$ and independent paths $R_2, R_2'$ in $G$ from $z_2$ 
to $D$ and internally disjoint from $V(p(K,{\cal A}))\cup V(X)$, such that $R_2$ ends at $t_2$, $R_2'$ ends at some vertex $t_2'\ne t_2$, and 
$y_1,t_2,t_2',w_2$ occur on $D$ in clockwise order. 

We may assume that 
\begin{itemize}
\item [(2)]  $K-\{u_1,u_2\}$ has no 2-separation $(K',K'')$ such that $V(K'\cap K'')\subseteq V(t_1Dt_2)$, $|V(K')|\ge 3$,
and $V(t_2Dt_1)\subseteq V(K'')$. 
\end{itemize}
For, suppose such a separation $(K',K'')$ does exist in $K-\{u_1,u_2\}$. Then by the definition of $u_1,u_2$, we see that $G$  has a separation $(G_1,G_2)$ such that $V(G_1\cap G_2)= V(K'\cap K'')\cup \{x_1, x_2, y_2\}$, $K'\subseteq V(G_1)$ and $K''\cup X\subseteq G_2$. Note that $G[\{x_1,x_2,y_2\}]$ is a triangle in $G$, 
$|V(G_2)|\ge 7$, and $|V(G_1)|\ge 6$ (as $|V(K')|\ge 3$). If $|V(G_1)|\ge 7$ then by Lemma~\ref{5cut_triangle}, $(i)$ or $(ii)$ or $(iii)$ holds. (Note that if $(iv)$ of 
Lemma~\ref{5cut_triangle} holds then $G'$ has a $TK_5$; so $(i)$ holds.)
So assume $|V(G_1)|=6$, and let $v\in V(G_1-G_2)$. Since $G$ is 5-connected, $N(v)=V(G_1\cap G_2)$. In particular, $v\ne y_1$ as $y_1y_2\notin E(G)$. 
Then $G[\{v,x_1,x_2,y_1\}]$ contains $K_4^-$, and $(ii)$ holds. So we may assume (2). 
\medskip

Next we may assume that 
\begin{itemize}
\item [(3)]  each neighbor of $x_1$ is contained in $V(X)$, or $V(t_1Dy_1)$, or some $A\in {\cal A}$ with $u_1\in N_K(A)$, and 
each neighbor of $x_2$ is contained $V(X)$, or $V(y_1Dt_2)$, or some $A\in {\cal A}$ with $u_2\in N_K(A)$. 
\end{itemize}
For, otherwise, we may assume by symmetry that there exists  $a\in N(x_1)-V(X)$ such that $a\notin V(t_1Dy_1)$ and
$a\notin A$ for $A\in {\cal A}$ with $u_1\in N_K(A)$.  
Let $a'=a$ and $S=a$ if $a\notin A$ for all $A\in {\cal A}$. When $a\in A$ for some $A\in {\cal  A}$ then by (2), there exists $a'\in N_K(A)-V(t_1Dt_2)$ 
and let $S$ be a path in $G[A +a']$ from $a$ to $a'$. By (2) again, there  
is a path $T$ from $a'$ to some $u\in V(t_2Dt_1)-\{t_1, t_2\}$ in $p(K,{\cal A})-\{u_1,u_2,y_2\}-t_1Dt_2$. 
Then $t_1Dt_2\cup R_1\cup R_2$ and $R_2'\cup t_2'Du\cup T$ give independent paths 
$T_1,T_2,T_3$ in $G-(X-\{z_1,z_2\})$ with $T_1, T_2$ from $y_1$ to $z_1, z_2$, respectively, and $T_{3}$ from $a'$ to $z_2$.  
Hence,  $z_2Xx_2\cup z_2Xy_2\cup T_2\cup (T_{3}\cup S\cup ax_1)\cup 
(T_1\cup z_1Xy_2)\cup G[\{x_1,x_2,y_1,y_2\}]$ is a $TK_5$ in $G'$ with branch vertices $x_1,x_2,y_1,y_2,z_2$; so $(i)$ holds.

\medskip

Label the vertices of $w_2Dy_1$ and $x_1Xy_2$ such that $w_2Dy_1=v_1\ldots v_k$ and $x_1Xy_2=v_{k+1}\ldots v_n$, with $v_1=w_2$, $v_k=y_1$, $v_{k+1}=x_1$ and $v_n=y_2$. 
Let $G_1$ denote the union of $x_1Xy_2$, $\{v_1, \ldots, v_k\}$, $G[A\cup (N_K(A)-u_1)]$ for $A\in {\cal A}$ with $u_1\in N_K(A)$,  all edges of $G'$
from $x_1Xy_2$ to $\{v_1,\ldots, v_k\}$, and all edges of $G'$ from $x_1Xy_2$ to $A$ for $A\in {\cal A}$ with $u_1\in N_K(A)$. 
Note that $G_1$ is $(4, \{v_1, \ldots, v_n\})$-connected. Similarly,  
let $y_1Dw_2=z_1\ldots z_l$ and $x_2Xy_2=z_{l+1}\ldots z_m$, with $z_1=w_2$, $z_l=y_1$, $z_{l+1}=x_2$ and $z_m=y_2$. 
Let $G_2$ denote the union of $y_2Xx_2$, $\{z_1, \ldots, z_l\}$, $G[A\cup (N_K(A)-u_2)]$ for $A\in {\cal A}$ with $u_2\in N_K(A)$,  all edges of $G'$
from $y_2Xx_2$ to $\{z_1,\ldots, z_l\}$, and all edges of $G'$ from $y_2Xx_2$ to $A$ for $A\in {\cal A}$ with $u_2\in N_K(A)$. 
Note that $G_2$ is $(4, \{z_1, \ldots, z_m\})$-connected. 

If both $(G_1,v_1,\ldots, v_n)$ and $(G_2, z_1, \ldots, z_m)$ are planar then $G-y_2$ is planar; so $(i)$ or $(ii)$ or $(iii)$ holds by Lemma~\ref{apex}. 
Hence, we may assume by symmetry that 
$(G_1,v_1,\ldots, v_n)$ is not planar. Then by Lemma~\ref{society}, there exist $1\le q<r<s<t\le n$ such that 
$G_1$ has disjoint paths $Q_1,Q_2$ from $v_q,v_r$ to $v_s,v_t$, respectively, and internally disjoint from $\{v_1, \ldots, v_n\}$. 

Since $(K,u_1,y_1,u_2,w_2)$ is 3-planar, it follows from the definition of $G_1$ that  
$q,r\le k$ and $s,t\ge k+1$. Note that the paths $y_1Dt_2$, $t_2'Dv_q$, $v_rDy_1$ give rise to independent paths $P_1,P_2,P_3$ in $K-\{u_1,u_2\}$, with 
$P_1$ from $y_1$ to $t_2$, $P_2$ from $t_2'$ to $v_q$, and $P_3$ from $v_r$ to $y_1$. 
Therefore,  $z_2Xx_2\cup z_2Xy_2\cup (R_2\cup P_1)\cup 
(R_2' \cup P_2\cup Q_1\cup v_sXx_1)\cup (P_3\cup Q_2 \cup v_tXy_2) \cup G[\{x_1,x_2,y_1,y_2\}]$ is a $TK_5$ in $G'$ with branch vertices $x_1,x_2,y_1,y_2,z_2$. 
So $(i)$ holds. \qed

\medskip

Conclusion $(iv)$ of Lemma~\ref{YZ} motivates the concept of $11$-tuple.
We say that $(G, X,x_1,x_2,$ $y_1,y_2,w_1,w_2,w_3,z_1,z_2)$ is an {\it $11$-tuple} if 
\begin{itemize}
\item $(G,X,x_1,x_2,y_1,y_2, w_1,w_2,w_3)$ is a $9$-tuple, and $z_i\in  V(x_iXy_2)-\{x_i,y_2\}$ for $i\in [2]$,
\item $H:=G'-(V(X-\{y_2,z_1,z_2\})\cup E(X))$ contains disjoint paths $Y,Z$ 
from $y_1,z_1$ to $y_2,z_2$, respectively, and
\item subject to the above conditions, $z_1Xz_2$ is maximal.
\end{itemize}

Since $G$ is 5-connected and $X$ is induced in $G'-x_1x_2$, each $z_i$ ($i\in [2]$) has
at least two neighbors in $H-\{y_2,z_1,z_2\}$ (which is 2-connected). Note that 
$y_2$ has exactly one neighbor $H-\{y_2,z_1,z_2\}$, namely, $w_2$. So $H-y_2$ is 2-connected. 
 
\begin{lem}\label{ABC}
Let $(G,X,x_1,x_2,y_1,y_2,w_1,w_2,w_3,z_1,z_2)$ be an $11$-tuple and    $Y,Z$ be disjoint paths in  $H:=G'-(V(X-\{y_2,z_1,z_2\})\cup E(X))$ from 
$y_1,z_1$ to $y_2,z_2$, respectively.
Then  $G$ contains \textsl{a} $TK_5$ in which $y_2$ is not a branch vertex, or $G'$ contains $TK_5$, or 
\begin{itemize}
 \item [$(i)$] for $i\in [2]$, $H$ has no path through $z_i,z_{3-i},y_1,y_2$ in order (so $y_1z_i\notin E(G)$), and 
\item [$(ii)$] there exists $i\in [2]$ such that $H$ contains independent paths $A,B,C$, with $A$ and $C$ from $z_i$ to $y_1$, and $B$  from $y_2$ to $z_{3-i}$.
\end{itemize}
\end{lem}
\pf First, suppose, for some $i\in [2]$, there is a path $P$ in $H$ from $z_i$  to $y_2$ such that $z_i,z_{3-i},y_1,y_2$ occur on $P$ in order. Then
$z_{3-i}Xx_{3-i}\cup z_{3-i}Xy_2\cup (z_{3-i}Pz_i\cup z_iXx_i)\cup z_{3-i}Py_1\cup y_1Py_2 \cup G[\{x_1,x_2,y_1,y_2\}]$ is a 
$TK_5$ with branch vertices $x_1,x_2,y_1,y_2,z_{3-i}$. So we
may assume that such $P$ does not exist. Hence by the existence of $Y,Z$ in
$H$, we have $y_1z_1,y_1z_2\notin E(G)$, and $(i)$ holds.

So from now on we may assume that $(i)$ holds. For each $i\in [2]$, let $H_i$ denote the graph obtained from
$H$ by duplicating $z_i$ and $y_1$, and let $z_i'$ and $y_1'$ denote the duplicates of $z_i$ and $y_1$, 
respectively. So in $H_i$, $y_1$ and $y_1'$ are not adjacent, and have the same set of neighbors, namely $N_H(y_1)$; and the same holds for $z_i$ and $z_i'$.

First, suppose for some $i\in [2]$,  $H_i$ contains pairwise disjoint paths $A',B',C'$ from $\{z_i,z_i',y_2\}$ to $\{y_1,y_1',z_{3-i}\}$, with
$z_i\in V(A'), z_i'\in V(C')$ and $y_2\in V(B')$. If $z_{3-i}\notin V(B')$, then after identifying $y_1$ with $y_1'$ and
$z_i$ with $z_i'$, we obtain from $A'\cup B'\cup C'$  a path in $H$
from $z_{3-i}$ to $y_2$ through $z_{i},y_1$ in order, contradicting
our assumption that $(i)$ holds. Hence $z_{3-i}\in V(B')$. Then we get the desired paths for $(ii)$
from $A'\cup B'\cup C'$ by identifying $y_1$ with $y_1'$ and $z_i$ with $z_i'$.

So we may assume that for each $i\in [2]$,
$H_i$ does not contain three pairwise disjoint paths from $\{y_2,z_i,z_i'\}$ to $\{y_1,y_1',z_{3-i}\}$.
Then $H_i$ has a separation $(H_i',H_i'')$ such that $|V(H_i'\cap H_i'')|=2$,
$\{y_2,z_i,z_i'\}\subseteq V(H_i')$ and $\{y_1,y_1',z_{3-i}\}\subseteq V(H_i'')$. 

We claim that $y_1,y_2,y_1',z_i',z_1,z_2\notin V(H_i'\cap H_i'')$ for $i\in [2]$. 
Note that $\{y_1,y_1'\}\ne V(H_i'\cap H_i'')$, since otherwise
$y_1$ would be a cut vertex in $H$ separating $z_{3-i}$ from
$\{y_2,z_i\}$. Now suppose one of $y_1,y_1'$ is in $V(H_i'\cap
H_i'')$; then since $y_1,y_1'$ are duplicates,
the  vertex in $V(H_i'\cap H_i'') -\{y_1,y_1'\} $ is a cut vertex in $H$
separating $\{y_1,z_{3-i}\}$ from $\{y_2,z_i\}$, a contradiction. So
$y_1,y_1'\notin V(H_i'\cap H_i'')$. Similar argument shows that
$z_i,z_i'\notin  V(H_i'\cap H_i'')$. Since $H-y_2$ is
2-connected, $y_2\notin  V(H_i'\cap H_i'')$. Since $H-\{z_{3-i},y_2\}$ is 2-connected, $z_{3-i}\notin 
V(H_i'\cap H_i'')$.

For $i\in [2]$, let $V(H_i'\cap H_i'')=\{s_i,t_i\}$, and let $F_i'$ (respectively, $F_i''$) be 
obtained from $H_i'$ (respectively, $H_i''$) by identifying $z_i'$ with
$z_i$ (respectively, $y_1'$ with $y_1$).  
Then $(F_i',F_i'')$ is a 2-separation in $H$ such that $V(F_i'\cap F_i'')=\{s_i,t_i\}$, 
$\{y_2,z_i\}\subseteq V(F_i')-\{s_i,t_i\}$, and $\{y_1,z_{3-i}\}\subseteq
V(F_i'')-\{s_i,t_i\}$. Let $Z_1,Y_2$ denote the $\{s_1,t_1\}$-bridges of 
$F_1'$ containing $z_1,y_2$, respectively; and let $Z_2,Y_1$ denote the
$\{s_1,t_1\}$-bridges of $F_1''$ containing 
$z_2,y_1$, respectively.

We may assume   $Y_1=Z_2$ or $Y_2=Z_1$. For, suppose $Y_1\ne Z_2$ and $Y_2\ne Z_1$. 
Since $H-y_2$ is 2-connected, there exist independent $P_1, Q_1$ in $Z_1$ from $z_1$ to $s_1, t_1$, respectively, independent 
paths $P_2, Q_2$ in $Z_2$ from $z_2$ to $s_1, t_1$, respectively,  independent paths $P_3, Q_3$ in $Y_1$ from $y_1$ to $s_1, t_1$, respectively, 
and a path $S$ in $Y_2$ from $y_2$ to one of  $\{s_1, t_1\}$ and avoiding the other, say avoiding $t_1$. 
Then $z_1Xx_1\cup z_1Xy_2\cup y_2x_1\cup  P_1\cup S\cup (P_3\cup y_1x_1)\cup (Q_2\cup Q_1)\cup P_2\cup  z_2Xy_2\cup (z_2Xx_2\cup x_2x_1)$
is a $TK_5$ in $G'$ with branch vertices $s_1,x_1,y_2,z_1,z_2$.  

Indeed, $Y_1=Z_2$. For, if $Y_1\ne Z_2$ then $Y_2=Z_1$, $Y_2-\{s_1,t_1\}$ has a path from $y_2$ to $z_1$, 
and $Y_1\cup Z_2$ has two independent paths from $y_1$ to $z_2$ (since $H-y_2$ is 2-connected). Now these three paths 
contradict the existence of the cut $\{s_2,t_2\}$ in $H$.  

Then $\{s_2,t_2\}\cap V(Y_1-\{s_1,t_1\})\ne \emptyset$. Without loss of generality, we may assume that $t_2\in V(Y_1)-\{s_1, t_1\}$. 
Suppose  $Y_2= Z_1$. Then $s_2\in V(Y_2)-\{s_1, t_1\}$ and we may assume that in $H$, $\{s_2,t_2\}$ separates $\{s_1, y_1, z_1\}$ from $\{t_1,y_2,z_2\}$. 
Hence, in $Y_1$, $t_2$ separates $\{y_1,s_1\}$ from $\{z_2,t_1\}$, and in $Y_2$, $s_2$ separates $\{z_1,s_1\}$ from $\{y_2,t_1\}$. 
But this contradicts the existence of the paths $Y$ and $Z$ in $H$.  So $Y_2\ne Z_1$. 
Since $H-y_2$ is 2-connected and $N_{G'}(y_2)=\{w_1,w_2,w_3,x_1,x_2\}$, we must have $s_2=w_2\in \{s_1, t_1\}$. By symmetry, we may assume that $s_2=w_2=s_1$. 

Let $Y_1', Z_2'$ be the $\{s_2,t_2\}$-bridge of $Y_1$ containing $y_1$, $z_2$, respectively. Then $t_1\notin V(Z_2')$; for, otherwise, $H-\{s_2,t_2\}$ 
would contain a path from $z_2$ to $z_1$, a contradiction. Therefore, 
because of the paths $Y$ and $Z$, $t_1\in V(Y_1')$ and $Y_1'$ contains disjoint paths $R_1,R_2$ from $s_2=s_1,t_1$ to $y_1,t_2$, respectively. Since $H-y_2$ is 2-connected, 
$Z_1$ has independent $P_1, Q_1$ from $z_1$ to $s_2=s_1, t_1$, respectively, and $Z_2'$ has independent 
paths $P_2, Q_2$ from $z_2$ to $s_2=s_1, t_2$, respectively. Now $z_1Xx_1\cup z_1Xy_2\cup y_2x_1\cup  P_1\cup s_1y_2\cup (R_1\cup y_1x_1)\cup P_2\cup (Q_2\cup R_2\cup Q_1)\cup 
z_2Xy_2\cup  (z_2Xx_2\cup x_2x_1)$
is a $TK_5$ in $G'$ with branch vertices $s_1,x_1,y_2,z_1,z_2$. \qed



\begin{lem}\label{PQ}
Let $(G,X,x_1,x_2,y_1,y_2,w_1,w_2,w_3,z_1,z_2)$ be an $11$-tuple and $Y,Z$ be disjoint paths in $H:=G'-V(X-\{y_2,z_1,z_2\} \cup E(X))$ from $y_1,z_1$ to $y_2,z_2$, respectively.  Then 
$G$ contains \textsl{a} $TK_5$ in which $y_2$ is not a branch vertex or $G'$ contains $TK_5$, or 
\begin{itemize}
 \item [$(i)$]  there exist $i\in [2]$ and independent paths $A,B,C$ in $H$, with $A$ and $C$ from $z_i$ to $y_1$, and $B$ from $y_2$ to $z_{3-i}$, 
\item [$(ii)$] for each $i\in [2]$ satisfying $(i)$, $z_{3-i}x_{3-i}\in E(X)$, and 
\item [$(iii)$] $H$ contains two disjoint paths from $V(B-y_2)$ to $V(A\cup C)-\{y_1,z_i\}$ and 
internally disjoint from $A\cup B\cup C$, with one ending in $A$ and the other ending in $C$.

\end{itemize}
\end{lem}

\pf By Lemma~\ref{ABC}, we may assume that 
\begin{itemize}
\item [(1)] for each $i\in [2]$, $H$ has no path through $z_i,z_{3-i},y_1,y_2$ in order (so $y_1z_i\notin E(G)$), and
\item [(2)] there exist $i\in [2]$ and independent paths $A,B,C$ in $H$, with $A$ and $C$ from $z_i$ to $y_1$, and $B$ from $y_2$ to $z_{3-i}$.
\end{itemize}

Let $J(A,C)$ denote the $(A\cup C)$-bridge of $H$ containing $B$, and $L(A,C)$ denote the union of
$(A\cup C)$-bridges of $H$ each of which intersects both $A-\{y_1,z_i\}$ and $C-\{y_1,z_i\}$. 
We choose $A, B, C$ such that the following are satisfied in the order listed:
\begin{itemize}
\item [(a)] $A, B,C$ are induced paths in $H$,
\item [(b)] whenever possible, $J(A,C)\subseteq L(A,C)$,
\item [(c)] $J(A, C)$ is maximal, and
\item [(d)] $L(A,C)$ is maximal.
\end{itemize}
We now show that $(ii)$ and $(iii)$ hold even with the restrictions (a), (b), (c) and (d) above. 
Let  $B'$ denote  the union of $B$ and the $B$-bridges of $H$ not
containing $A\cup C$.

\begin{itemize}
\item [(3)] If $(iii)$ holds then $(ii)$ holds.
\end{itemize}
Suppose $(iii)$ holds. Let $V(P \cap B)=\{p\}$, $V(Q \cap B)=\{q\}$, $V(P \cap C)=\{c\}$ and $V(Q \cap A)=\{a\}$.  
By the symmetry between $A$ and $C$, we may assume that $y_2,p,q,z_{3-i}$ occur on $B$ in order. 
We may further choose $P,Q$ so that $pBz_{3-i}$ is maximal.

To prove $(ii)$, suppose there exists $x\in V(z_{3-i}Xx_{3-i})-\{x_{3-i},z_{3-i}\}$. If $N(x)\cap V(H)-\{y_1\}\not\subseteq V(B')$ then $G'$ has a path $T$ from $x$ to 
$(A-y_1)\cup (C-y_1)\cup (P-p)\cup (Q-a)$ and internally disjoint from $A\cup B'\cup C\cup P\cup Q$; so $A\cup B\cup C\cup P\cup Q\cup T$  contain disjoint paths 
from $y_1,z_i$ to $y_2,x$, respectively, contradicting the choice of $Y$ and $Z$ in the $11$-tuple (that $z_1Xz_2$ is maximal). 
So  $N(x)\cap V(H)-\{y_1\}\subseteq V(B')$. Consider $B'':=G[(B'-z_{3-i})+x]$. 

If $B''$ contains disjoint paths $P',Q'$ from 
$y_2,x$ to $p,q$, respectively, then $Q'\cup Q\cup
aAz_i$ and  $P'\cup P\cup cCy_1$ contradict the choice of $Y,Z$.
If $B''$ contains disjoint paths $P'',Q''$ from $x,y_2$ to $p,q$, respectively, then 
 $Q''\cup Q\cup aAy_1$ and $P''\cup P\cup cCz_i$ contradict the choice of
$Y,Z$. 

So we may assume that there is a cut vertex
$z$ in $B''$ separating $\{x,y_2\}$ from $\{p,q\}$. Note that $z\in V(y_2Bp)$.

Since $x$ has at least two neighbors in $B''-y_2$ (because $G$ is $5$-connected
and $X$ is induced in $G'-x_1x_2$), the $z$-bridge of $B''$ containing $\{x,y_2\}$ has at least three vertices. 
 Therefore, from the maximality of $pBz_{3-i}$ and 2-connectedness of $H-\{y_2,z_1,z_2\}$, there is a path in
$H$ from $y_1$ to $y_2Bz-\{y_2,z\}$ and internally disjoint 
from $P\cup Q\cup A\cup C\cup B'$. So there is a path
$Y'$ in $H$
from $y_1$ to $y_2$ and disjoint from  $P\cup Q\cup A\cup C\cup pBz_{3-i}$. Now 
$z_{3-i}Bp\cup P\cup cCz_i\cup A\cup Y'$  is a path in $H$ through
$z_{3-i},z_i,y_1,y_2$ in order, contradicting (1).  

\medskip

By (2) and (3), it suffices to prove $(iii)$.  
Since $H-\{y_2,z_i\}$ is 2-connected, it contains disjoint paths
$P,Q$ from $B-y_2$ to some distinct vertices $s,t\in V(A\cup C)-\{z_i\}$, respectively, and internally
disjoint from $A\cup B\cup C$.

\begin{itemize}
\item [(4)] We may choose $P,Q$ so that  $s\ne y_1$ and $t\ne y_1$. 
\end{itemize}
 For, otherwise, $H-\{y_2,z_i\}$ has a separation $(H_1,H_2)$ such that
$V(H_1\cap H_2)=\{v,y_1\}$ for some $v\in V(H)$, $(A\cup
C)-z_i\subseteq H_1$ and $B-y_2\subseteq H_2$. 
Recall the disjoint paths
$Y,Z$ in $H$ from $z_1,y_1$ to $z_2,y_2$, respectively. 
Suppose $v\notin V(Z)$. Then $Z-z_i\subseteq H_2-\{y_1,v\}$. Hence we may choose  $Y$ (by modifying $Y\cap H_1$) so that $V(Y\cap A)=\{y_1\}$ or $V(Y\cap
C)=\{y_1\}$. Now $Z\cup A\cup Y$ or $Z\cup C\cup Y$ is a path in $H$ from $z_{3-i}$ to $y_2$ through $z_i,y_1$ in order,
contradicting (1). So $v\in V(Z)$. 
 Hence $Y\subseteq H_2-v$, and we may choose  $Z$ (by modifying $Z\cap H_1$) so that $V(Z\cap A)=\{z_i\}$ or $V(Z\cap
C)=\{z_i\}$. Now $Z\cup A\cup Y$ or $Z\cup C\cup Y$ is a path in $H$ from $z_{3-i}$ to $y_2$ through $z_i,y_1$ in order, 
contradicting (1)  and completing the proof of (4).

\medskip

If $s\in V(A-y_1)$ and $t\in V(C-y_1)$ or $s\in V(C-y_1)$ and $t\in V(A-y_1)$,
then $P,Q$ are the desired paths for $(iii)$. So we may assume by
symmetry that $s,t\in V(C)$. Let
$V(P\cap B)=\{p\}$ and $V(Q\cap B)=\{q\}$ such that $y_2,p,q,z_{3-i}$ occur on $B$ in this order.
 By (1) $z_i,s,t,y_1$ must occur on $C$ in order. We choose $P,Q$ so that 

\begin{itemize}
\item [$(*)$] $sCt$ is maximal, then $pBz_{3-i}$ is maximal, and then $qBz_{3-i}$ is minimal.
\end{itemize} 

Now consider $B'$, the union of $B$ and the $B$-bridges of $H$ not containing
$A\cup C$. Note that $(P-p)\cup (Q-q)$ is disjoint from $B'$, and every
path in $H$ from $A\cup C$ to $B'$ and internally disjoint from $A\cup B'\cup
C$ must end in $B$. For convenience, let $K=P\cup Q\cup A\cup B'\cup C$.

\begin{itemize}
\item [(5)]  $B'-y_2$ contains independent paths $P',Q'$ from  $z_{3-i}$ to $p,q$, respectively.
\end{itemize}
Otherwise, $B'-y_2$ has a cut vertex $z$ separating $z_{3-i}$ 
from $\{p,q\}$. Clearly, $z\in V(qBz_{3-i}-z_{3-i})$, and we choose $z$ so that
$zBz_{3-i}$ is minimal. 

Let $B''$ denote the $z$-bridge of $B'-y_2$ containing $z_{3-i}$; then $zBz_{3-i}\subseteq B''$.  
Since $H-\{y_2,z_i\}$ is
2-connected, it  contains a path  $W$ from some $w'\in V(B''-z)$ to some 
$w\in V(P\cup Q\cup A\cup C)-\{z_i\}$ and internally disjoint 
from $K$. By the definition of $B'$, $w'\in V(z_iBz_{3-i})$. By (1), $w\notin V(P)\cup V(z_iCt-t)$. By  $(*)$, $w\notin V(Q)\cup V(tCy_1-y_1)$. 

If $w\in V(A)-\{z_i,y_1\}$ then $P,W$ give the desired paths for $(iii)$. So we may assume $w=y_1$ for any choice of $W$; hence,
$z\in V(Z)$ and $Y\cap (B''\cup (W-y_1))=\emptyset$.  
 By the minimality of $zBz_{3-i}$, $B''$ has independent paths  
$P'',Q''$ from $z_{3-i}$ to $z,w'$, respectively. Note that $z_iZz\cap (B''-z)=\emptyset$. 
Now $z_iZz\cup P''\cup Q''\cup W\cup Y$  
is a path in $H$  through $z_i,z_{3-i},y_1,y_2$ in order, contradicting (1).

\begin{itemize}
\item [(6)] We may assume that $J(A,C)\not\subseteq L(A,C)$. 
\end{itemize}
 For, otherwise, there is a path
$R$ from $B$ to some $r\in V(A)-\{y_1,z_i\}$ and internally disjoint from $A\cup B'\cup C$.
If $R\cap (P\cup Q)\ne \emptyset$, then it is easy to check that $P\cup Q\cup R$
contains the desired paths for $(iii)$. So we may assume $R\cap (P\cup
Q)=\emptyset$. If $y_2\notin V(R)$, then $P,R$ are the desired paths for $(iii)$. So assume
$y_2\in V(R)$. Recall the paths $P',Q'$ from (5). Then  $z_iCs\cup P\cup P'\cup Q'\cup Q\cup tCy_1\cup y_1Ar\cup R$ is a 
path in $H$  through $z_i,z_{3-i},y_1,y_2$ in order, contradicting (1) and completing the proof of (6). 
 
\medskip

Let $J=J(A,C)\cup C$. Then by (1), $J$ does not contain disjoint paths from $y_2,z_i$ to $y_1,z_{3-i}$, respectively. So 
by Lemma~\ref{2path}, there exists a
collection ${\cal A}$ of subsets of $V(J)-\{y_1,y_2,z_1,z_2\}$ such
that $(J,{\cal A}, z_i,y_1,z_{3-i},y_2)$ is 3-planar.  We choose ${\cal A}$ so that every member of ${\cal A}$ is minimal and, subject to this, 
$|{\cal A}|$ is minimum. Then

\begin{itemize}
\item [(7)]  for any $D\in {\cal A}$ and any $v\in V(D)$,  $(J[D+N_J(D)],N_J(D)\cup \{v\})$ is not 3-planar. 
\end{itemize}
Suppose for some $D\in {\cal A}$ and some $v\in D$, there is a collection of subsets
${\cal A}'$ of $D-\{v\}$ such that $(J[D+N_J(D)],{\cal A}',N_J(D)\cup \{v\})$ is 3-planar. Then, with ${\cal A}''=({\cal A}-\{D\})\cup 
{\cal A}'$, $(J, {\cal A}'', z_i,y_1,z_{3-i},y_2)$ is 3-planar. So ${\cal A}''$ contradicts the choice of ${\cal A}$. 
Hence, we have (7). 

\medskip

 Let $v_1,\ldots,v_k$ be the vertices of $L(A,C)\cap (C-\{y_1,z_i\})$ such that $z_i,v_1,\ldots, v_k,y_1$ 
occur on $C$ in the order listed. We claim that 

\begin{itemize}
\item [(8)] $(J,z_i,v_1,\ldots, v_k,y_1,z_{3-i},y_2)$ is 3-planar.
\end{itemize}
For, suppose otherwise.  Since there is only one
$C$-bridge in $J$ and $(J,{\cal A},z_i,y_1,z_{3-i},y_2)$ is
3-planar,  there exist $j\in [k]$ and $D\in {\cal A}$ such that $v_j\in D$.  Since $H$ is 2-connected,  
let  $c_1,c_2\in V(C)\cap N_J(D)$ with  $c_1Cc_2$ maximal. 

Suppose $N_J(D)\subseteq V(C)$. Then, since there is only one
$C$-bridge in $J$ and $(J,{\cal A},z_i,y_1,z_{3-i},y_2)$ is
3-planar, $J$ has a separation $(J_1,J_2)$ such that $V(J_1\cap J_2)=\{c_1,c_2\}$, $D\cup V(c_1Cc_2)\subseteq V(J_1)$, and $B\subseteq J_2$. 
Since $J$ has only one $C$-bridge and $C$ is induced in $H$, we have $J_1=c_1Cc_2$. Now let ${\cal A}'$ be obtained from ${\cal A}$ by removing all 
members of ${\cal A}$ contained in $V(J_1)$. Then $(J, {\cal A}', z_i, y_1,z_{3-i},y_2)$ is 3-planar, contradicting the choice of ${\cal A}$. 
 
Thus, let $c\in N_J(D)- V(C)$. So $c\in V(J(A,C))$. Let $D'=J[D+ \{c_1,c_2,c\}]$.
By (7) and Lemma~\ref{2path}, $D'$ contains disjoint paths $R$
from $v_j$ to $c$ and $T$ from $c_1$ to $c_2$. We may assume $T$ is
induced. Let $C'$ be obtained from $C$ by replacing $c_1Cc_2$ with
$T$. We now see that $A,B,C'$ satisfy (a), but $J(A, C')$ intersects both $A-\{y_1,z_i\}$ (by definition of $v_j$ and because $c\in V(J(A,C))-V(C)$) and $C'-\{y_1,z_i\}$
(because of $P,Q$), contradicting (b) (via (6)) and completing the
proof of (8).

\begin{itemize}
\item [(9)] There exist disjoint paths $R_1, R_2$ in $L(A,C)$ from some $r_1,r_2\in
V(C)$ to some $r_1',r_2'\in V(A)$, respectively, and internally disjoint from $A\cup C$, such that
$z_i,r_1,r_2,y_1$ occur on $C$ in this order and $z_i,r_2',r_1',y_1$ occur on $A$ in this order.
\end{itemize}
We prove (9) by studying the $(A\cup C)$-bridges of $H$ other than $J(A,C)$. 
For any $(A\cup C)$-bridge $T$ of $H$ with $T\ne J(A,C)$, if $T$ intersects $A$ let $a_1(T), a_2(T)\in V(T\cap A)$ with 
$a_1(T)Aa_2(T)$ maximal, and if $T$ intersects  $C$ let $c_1(T), c_2(T)\in V(T\cap C)$ with $c_1(T)Cc_2(T)$ maximal. We choose the notation so
that $z_i,a_1(T),a_2(T), y_1$ occur on $A$ in order, and
$z_i,c_1(T),c_2(T), y_1$ occur on $C$ in order.

 If $T_1,T_2$ are $(A\cup
C)$-bridges of $H$  such that $T_2\subseteq L(A,C)$,  $T_1\ne J(A,C)$, and $T_1$
intersects $C$ (or $A$) only, then
$c_1(T_1)Cc_2(T_1)-\{c_1(T_1),c_2(T_1)\}$ (or
$a_1(T_1)Aa_2(T_1)-\{a_1(T_1),a_2(T_1)\}$) does not intersect 
$T_2$.  For, otherwise, we may modify $C$ (or $A$) by replacing
$c_1(T_1)Cc_2(T_1)$ (or  $a_1(T_1)Aa_2(T_1)$) with an induced path
in $T_1$ from $c_1(T_1)$ to $c_2(T_1)$ (or from $a_1(T_1)$ to
$a_2(T_1)$). The new $A$ and $C$ do not affect (a), (b) and (c) but enlarge $L(A,C)$, 
contradicting (d). 

Because of the disjoint paths $Y$ and $Z$ in $H$, $(H,z_i,y_1,z_{3-i},y_2)$ is not 3-planar. 
By (1) $A-\{y_1,z_i\}\ne \emptyset$. Hence, since $H-\{y_2,z_1,z_2\}$ is 2-connected, $L(A,C)\ne \emptyset$. 
Thus, since $(J,z_i,v_1,\ldots,$ $ v_k,y_1,z_{3-i},y_2)$ is 3-planar (by (8)) and $J(A,C)$ 
does not intersect $A-\{y_1,z_i\}$ (by (6)), one of the following holds: There
exist $(A\cup C)$-bridges $T_1,T_2$ of $H$ such
that $T_1\cup T_2\subseteq L(A,C)$, $z_iAa_2(T_1)$ properly contains
$z_iAa_1(T_2)$, and  $c_1(T_1)Cy_1$ properly
contains $c_2(T_2)Cy_1$; or there exists an $(A\cup C)$-bridge
$T$ of $H$ such that $T\subseteq L(A,C)$ and $T\cup a_1(T)Aa_2(T)\cup
 c_1(T)Cc_2(T)$ has disjoint paths from $a_1(T),a_2(T)$ to $c_2(T),c_1(T)$, respectively.
In either case, we have (9). 

\begin{itemize}
\item [(10)] $r_1,r_2\in V(tCy_1)$ for all choices of $R_1,R_2$ in (9),   or $r_1,r_2\in V(z_iCs)$  for all choices of $R_1,R_2$ in (9). 
\end{itemize}
For, suppose there exist $R_1,R_2$ such that $r_1\in V(z_iCs)$ and $r_2\in V(tCy_1)$, or
$r_1\in V(sCt)-\{s,t\}$, or $r_2\in V(sCt)-\{s,t\}$. Let $A':=z_iAr_2'\cup R_2\cup r_2Cy_1$ and
$C':=z_iCr_1\cup R_1\cup r_1'Ay_1$. We may assume $A',C'$ are induced paths
in $H$ (by taking  induced paths in $H[A']$ and $H[C']$). 
Note that $A',B,C'$ satisfy (a), and $J(A,C)\subseteq J(A',C')$. 
However, because of $P$ and $Q$,  $J(A',C')$ intersects both $A'-\{z_i,y_1\}$ and $C'-\{z_i,y_1\}$,
contradicting (b) (via (6)) and completing the proof of (10). 

\medskip

If $r_1, r_2\in V(z_iCs)$ for all choices of  $R_1,R_2$ in (9) then we choose such $R_1,R_2$ that 
$z_iAr_1'$ and $z_iCr_2$ are maximal, and let $z':=r_1'$ and $z''=r_2$; otherwise, define $z'=z''=z_i$.
Similarly, if  $r_1, r_2\in V(tCy_1)$ for all choices of $R_1,R_2$ in (9), then we choose such $R_1,R_2$ that 
$y_1Ar_2'$ and $y_1Cr_1$ are maximal, and let $y':=r_2'$ and $y''=r_1$; otherwise, define $y'=y''=y_1$.
By (10), $z_i,z',y',y_1$ occur on $A$ in order, and $z_i,z'',s,t,y'',y_1$ occur on $C$ in order. 

Note that $H$ has a path $W$ from some $y\in V(B)\cup V(P-s)\cup V(Q-t)$ 
to some $w\in V(z_iAz'-\{z',z_i\})\cup V(z_iCz''-\{z'',z_i\})\cup V(y'Ay_1-\{y',y_1\})\cup V(y''Cy_1-\{y'',y_1\})$ such that  
$W$ is internally disjoint from $K$. For, otherwise, 
$(H,z_i,y_1,z_{3-i},y_2)$ is 3-planar, contradicting the existence of the disjoint paths $Y$ and $Z$. 
By (6),  $w\notin V(A)$. If $w\in V(z_iAz'-\{z',z_i\})\cup V(y'Ay_1-\{y',y_1\})$  then we can find the desired $P,Q$.
So assume   $w\in V(z_iCz''-\{z'',z_i\})\cup V(y''Cy_1-\{y'',y_1\})$. 
 By $(*)$ and (1), $y\notin V(B-y_2)$ and $y\notin V(P\cup Q)$. 
This forces $y=y_2$, which is impossible as $N_H(y_2)=\{w_2\}$. \qed

\medskip

{\it Remark}. Note from the proof of Lemma~\ref{PQ} that the conclusions $(ii)$ and $(iii)$ hold for those paths $A,B,C$ that satisfy (a), (b), (c) and (d). 

\section{Finding $TK_5$}

In this section, we prove Theorem~\ref{y_2}. Let $G$ be a 5-connected nonplanar graph and let $x_1,x_2,y_1,y_2\in V(G)$ be distinct such that 
$G[\{x_1,x_2,y_1,y_2\}]\cong K_4^-$ and $y_1y_2\notin E(G)$. Let $w_1,w_2,w_3\in N(y_2)-\{x_1,x_2\}$ be distinct and let  
$G':=G-\{y_2v:v\notin \{w_1,w_2,w_3,x_1,x_2\}\}.$ 

We may assume that $G'-x_1x_2$ has an induced path $L$ from $x_1$ to $x_2$ such that $y_1,y_2\notin V(L)$,  $(G-y_2)-L$ is 2-connected, 
and $w_1,w_2,w_3\in V(L)$; for otherwise, the conclusion of Theorem~\ref{y_2} follows from Lemma~\ref{reduction}. 
Hence, $G'-x_1x_2$ has an induced path $X$ from $x_1$ to $x_2$ such that $y_1\notin V(X)$, $w_1y_2,w_3y_2\in E(X)$, 
and $G'-X=G-X$ is 2-connected. Hence, $(G,X,x_1,x_2,y_1,y_2,w_1,w_2,w_3)$ is a 9-tuple. 

We may assume that there exist $z_i\in V(x_iXy_2)-\{x_i,y_2\}$ for $i\in [2]$ such that $H:=G'-(X-\{y_2,z_1,z_2\})$ has disjoint paths $Y,Z$ from 
$y_1,z_1$ to $y_2,z_2$, respectively; for, otherwise, the conclusion of Theorem~\ref{y_2} follows from  Lemma~\ref{YZ}. 
We choose such $Y,Z$ so that $z_1Xz_2$ is maximal. Then $(G,X,x_1,x_2,y_1,y_2,w_1,w_2,w_3,z_1,z_2)$ is an 11-tuple.

\begin{figure}
\begin{center}
\includegraphics[scale=0.25]{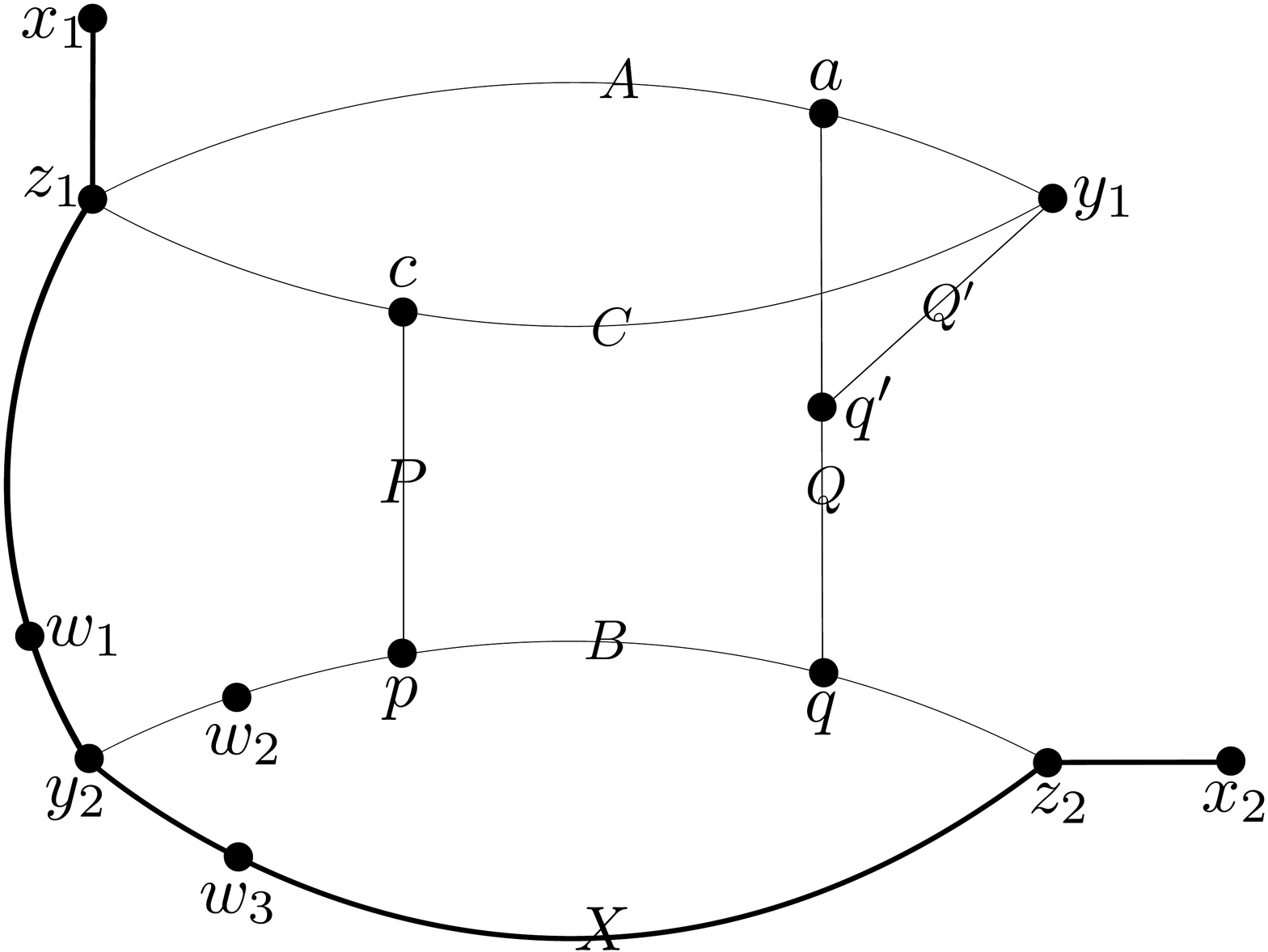}
\caption{\label{structure} An intermediate structure}
\end{center}
\end{figure}

 By Lemma~\ref{ABC} and by symmetry, we may assume that 
\begin{itemize}
\item [(1)] for $i\in [2]$, $H$ has no path through $z_i,z_{3-i},y_1,y_2$ in order (so $y_1z_i\notin E(G)$), 
\end{itemize} 
and that there exist independent paths $A,B,C$ in $H$ 
with $A$ and $C$ from $z_1$ to $y_1$, and $B$ from $y_2$ to
$z_2$.  See Figure~\ref{structure}.

Let $J(A,C)$ denote the $(A\cup C)$-bridge of $H$ containing $B$, and $L(A,C)$ denote the union of $(A\cup C)$-bridges of $H$ 
intersecting both $A-\{y_1,z_1\}$ and  $C-\{y_1,z_1\}$.
We may choose $A, B, C$ such that 
the following are satisfied in the order listed:
\begin{itemize}
\item [(a)] $A, B,C$ are induced paths in $H$,
\item [(b)] whenever possible $J(A, C)\subseteq L(A,C)$, 
\item [(c)] $J(A, C)$ is maximal, and 
\item [(d)] $L(A, C)$ is maximal.
\end{itemize}

By Lemma~\ref{PQ} and its proof (see the remark at the end of Section 4), we may assume that $$z_2x_2\in E(X)$$ 
and that there exist disjoint paths $P,Q$ in $H$ from $p,q\in V(B-y_2)$ to
  $c\in V(C)-\{y_1,z_1\}, a\in V(A)-\{y_1,z_1\}$,
respectively, and internally disjoint from $A\cup B\cup C$. 
By  symmetry between $A$ and $C$, we assume that
$y_2,p,q,z_2$ occur on $B$ in order. We further choose $A,B,C,P,Q$ so that

\begin{itemize}
\item [(2)] $qBz_2$ is minimal, then $pBz_2$ is maximal, and then  
$aAy_1\cup cCz_1$ is minimal.
\end{itemize}

Let $B'$ denote the union of $B$ and
the $B$-bridges of $H$ not containing $A\cup C$. Note that all paths in $H$ from $A\cup C$ to $B'$ and 
internally disjoint from $B'$ must have an end in $B$. 
For convenience, let $$K:= A\cup B'\cup C\cup P\cup Q.$$ 
Then

\begin{itemize}
\item [(3)] $H$ has no path from $aAy_1-a$ to $z_1Cc-c$ and internally disjoint from $K$. 
\end{itemize}
For, suppose $S$ is a path in $H$ from some vertex  $s\in V(aAy_1-a)$ to some vertex $s'\in
V(z_1Cc-c)$ and internally disjoint from $K$.
Then $z_2Bq\cup Q\cup aAz_1\cup z_1Cs'\cup S\cup sAy_1\cup
y_1Cc\cup P\cup pBy_2$ is a path in $H$ through $z_2,z_1,y_1,y_2$ in order,  contradicting (1).

\medskip

We proceed by proving a number of  claims from which Theorem~\ref{y_2} will
follow. Our intermediate goal is to prove (12) that $H$ contains a path from $y_1$ to $Q-a$ and internally disjoint from $K$.
However, the claims leading to (12) will also be useful when we later consider structure of $G$ near $z_1$. 
\medskip

\begin{itemize}
\item [(4)] $B'-y_2$ has no cut vertex contained in $qBz_2-z_2$ and, hence, 
for any $q^*\in V(B')-\{y_2,q\}$,  $B'-y_2$ has independent paths $P_1,P_2$ from $z_2$ to $q,q^*$, respectively.
\end{itemize}
Suppose $B'-y_2$ contains a cut vertex $u$ with  $u\in V(qBz_2-z_2)$. Choose $u$ so that 
$uBz_2$ is minimal. 
Since $H-\{y_2,z_1\}$ is 2-connected, there is a path $S$ in $H$ from some 
$s'\in V(uBz_2-u)$ to some $s\in  V(A\cup C\cup P\cup Q)-\{p,q\}$ and internally disjoint from
$K$. By the minimality of $uBz_2$, the 
$u$-bridge of  $B'-y_2$  containing $uBz_2$ has independent
paths $R_1,R_2$ from $z_2$ to $s',u$, respectively.
 By the minimality of $qBz_2$ in (2), $S$ is disjoint from $(P\cup Q \cup A \cup C) -\{z_1,y_1\} $. 
If $s=z_1$ then $(R_1\cup S)\cup A\cup (y_1Cc\cup P\cup pBy_2)$ is a path in $H$ through $z_2,z_1,y_1,y_2$ in order, contradicting  (1). 
So $s=y_1$.  Then $(z_1Aa\cup Q\cup qBu\cup R_2)\cup (R_1\cup S)\cup (y_1Cc\cup P\cup pBy_2)$ is a path in
$H$ through $z_1,z_2,y_1,y_2$ in order, contradicting (1).  

Hence, $B'-y_2$ has no cut vertex contained in $qBz_2-z_2$. Thus, the second half of (4) follows from Menger's theorem.

\begin{itemize} 
\item [(5)] We may assume that $G'$ has no path from $aAy_1-a$ to $z_1Xz_2$ and internally disjoint from $K\cup X$, and no path from $cCy_1-c$ 
to $z_1Xz_2-z_1$ and internally disjoint from $K\cup X$.
\end{itemize}
For, suppose $S$ is a path in $G'$ from some $s\in V(aAy_1-a)\cup V(cCy_1-c)$ to some $s'\in V(z_1Xz_2)$ and  internally disjoint from $K\cup X$, such that 
$s'\ne z_1$ if $s\in  V(cCy_1-c)$. If $s'=z_1$ then  $s\in V(aAy_1-a)$; so $z_2Bq\cup Q\cup aAz_1\cup S\cup sAy_1\cup y_1Cc\cup P\cup pBy_2$ is a 
path in $H$  through $z_2,z_1,y_1,y_2$ in order, contradicting (1). 
If $s'=z_2$ then $s=y_1$ by (2); so  $(z_1Aa\cup Q\cup qBz_2)\cup S\cup y_1Cc\cup P\cup pBy_2$ 
is a path in $H$  through $z_1,z_2,y_1,y_2$ in order, contradicting (1).   
Hence,  $s'\in V(z_1Xz_2)-\{z_1, z_2\}$.

Suppose  $s'\in V(z_1Xy_2-z_1)$. Let $P_1,P_2$ be the paths in (4) with $q^*=p$. If $s\in V(aAy_1-a)$ then $z_2x_2\cup z_2Xy_2\cup (P_2\cup P\cup cCy_1)\cup 
(P_1\cup Q\cup aAz_1\cup z_1Xx_1)\cup (y_1As\cup S\cup s'Xy_2)\cup 
G[\{x_1,x_2,y_1,y_2\}]$ is a $TK_5$ in $G'$ with branch vertices $x_1,x_2,y_1,y_2,z_2$.
If $s\in V(cAy_1-c)$ then $z_2x_2\cup z_2Xy_2\cup (P_2\cup P\cup cCz_1\cup z_1Xx_1)\cup 
(P_1\cup Q\cup aAy_1)\cup (y_1Cs\cup S\cup s'Xy_2)\cup 
G[\{x_1,x_2,y_1,y_2\}]$ is a $TK_5$ in $G'$ with branch vertices $x_1,x_2,y_1,y_2,z_2$.

Now assume $s'\in V(z_2Xy_2-z_2)$.
If $s\in V(aAy_1-a)$, then $z_1Xx_1\cup z_1Xy_2\cup C\cup (z_1Aa\cup Q\cup qBz_2\cup z_2x_2)\cup 
(y_1As\cup S\cup s'Xy_2)\cup 
G[\{x_1,x_2,y_1,y_2\}]$ is a $TK_5$ in $G'$ with branch vertices $x_1,x_2,y_1,y_2,z_1$.
If $s\in V(cCy_1-c)$, then $z_1Xx_1\cup z_1Xy_2\cup A\cup (z_1Cc\cup P\cup pBz_2\cup z_2x_2)\cup 
(y_1Cs\cup S\cup s'Xy_2)\cup 
G[\{x_1,x_2,y_1,y_2\}]$ is a $TK_5$ in $G'$ with branch vertices $x_1,x_2,y_1,y_2,z_1$. 
This completes the proof of (5).
\medskip

Denote by $L(A)$ (respectively, $L(C)$) the union of $(A\cup C)$-bridges of $H$ not intersecting $C$ (respectively, $A$). 
Let $C'= C\cup L(C)$. The next four claims concern paths from $x_1Xz_1-z_1$ to other parts of $G'$. 
 We may assume that

\begin{itemize}
\item [(6)]  $N(x_1Xz_1-\{x_1, z_1\})\subseteq V(C')\cup \{x_1, z_1\}$, and  that $G'$ has no disjoint paths from $s_1, s_2\in V(x_1Xz_1-z_1)$ 
to $s_1', s_2'\in V(C)$, respectively, and internally disjoint from $K\cup X$ such that $s_2'\in V(cCy_1-c)$, $x_1, s_1, s_2, z_1$ occur on $X$ in order, 
and $z_1, s_1', s_2', y_1$ occur on $C$ in order.
\end{itemize}
First,  suppose  $N(x_1Xz_1-\{x_1, z_1\})\not \subseteq V(C')\cup \{x_1, z_1\}$. Then there exists a path $S$ in $G'$ from some  $s\in
V(x_1Xz_1)-\{x_1, z_1\}$ to some $s'\in V(A\cup B'\cup P\cup Q)-\{c,y_1,y_2,z_1, z_2\}$ and internally
disjoint from $K\cup X$. If $s'\in V(A)-\{z_1, y_1\}$ then $y_1Cc\cup P\cup pBy_2$, $S\cup s'Aa\cup Q\cup qBz_2$ 
contradict the choice of $Y$, $Z$. If $s'\in V(Q-a)$ then $y_1Cc\cup P\cup pBy_2$, $S\cup s'Qq\cup qBz_2$ contradict the choice of $Y$, $Z$. 
If $s'\in V(P-c)$ then let $P_1,P_2$ be the paths in (4) with $q^*=p$; now 
$ z_2x_2\cup z_2Xy_2\cup (P_1\cup Q\cup aAy_1)\cup (P_2\cup pPs'\cup S\cup sXx_1)\cup (C\cup z_1Xy_2)
\cup G[\{x_1, x_2, y_1, y_2\}]$ is a $TK_5$ in $G'$ with branch vertices $x_1, x_2, y_1, y_2, z_2$. 
If $s'\in V(B')-\{y_2, p, q\}$ then let $P_1,P_2$ be the paths in (4)  with $q^*=s'$; now $z_2x_2\cup z_2Xy_2\cup (P_1\cup Q\cup aAy_1)\cup (P_2\cup S\cup sXx_1)
\cup (C\cup z_1Xy_2)\cup G[\{x_1, x_2, y_1, y_2\}]$ is a $TK_5$ in $G'$ with branch vertices 
$x_1, x_2, y_1, y_2, z_2$.

Now assume $G'$ has disjoint paths $S_1, S_2$ from $s_1, s_2\in V(x_1Xz_1-z_1)$ 
to $s_1', s_2'\in V(C)$, respectively, and internally disjoint from $K\cup X$ such that $s_2'\in V(cCy_1-c)$, $x_1, s_1, s_2, z_1$ occur on $X$ in order, 
and $z_1, s_1', s_2', y_1$ occur on $C$ in order. Let $P_1,P_2$ be the paths in (4)  with $q^*=p$. 
Then $z_2x_2\cup z_2Xy_2\cup (P_1\cup Q\cup aAy_1)\cup (P_2\cup P\cup cCs_1'\cup S_1\cup s_1Xx_1)
\cup (y_1Cs_2'\cup S_2\cup s_2Xy_2)\cup G[\{x_1, x_2, y_1, y_2\}]$ 
is a $TK_5$ in $G'$ with branch vertices $x_1, x_2, y_1, y_2, z_2$. 
This completes the proof of (6).

\begin{itemize}
\item [(7)] For any path $W$ in $G'$ from $x_1$ to some $w\in V(K)-\{y_1,z_1\}$ and internally disjoint from $K\cup X$, we may assume  
$w\in V(A\cup C)-\{y_1, z_1\}$. (Note that such $W$ exists as $G$ is 5-connected and $G'-X$ is 2-connected.)
\end{itemize}
For, let $W$ be a  path in $G'$ from $x_1$ to $w\in V(K)-\{y_1,z_1\}$ and internally disjoint from $K\cup X$, 
such that $w\notin V(A\cup C)-\{z_1, y_1\}$. Then $w\ne y_2$ as $N_{G'}(y_2)=\{w_1,w_2,w_3,x_1,x_2\}$.

Suppose $w\in V(B'-q)$. Let $P_1,P_2$ be the paths in (4) with $q^*=w$. Then  
$z_2x_2\cup z_2Xy_2\cup  (P_1\cup Q\cup aAy_1)\cup (P_2\cup W)\cup (C\cup z_1Xy_2)\cup G[\{x_1, x_2, y_1, y_2\}]$ is 
a $TK_5$ in $G'$ with branch vertices $x_1, x_2, y_1, y_2, z_2$.

So assume $w\notin V(B'-q)$. Let $P_1,P_2$ be the paths in (4) with $q^*=p$.
If  $w\in V(P-c)$ then $z_2x_2\cup z_2Xy_2\cup (P_1\cup Q\cup aAy_1)\cup (P_2\cup pPw\cup W)\cup (C\cup z_1Xy_2)\cup 
G[\{x_1,y_1,x_2,y_2\}]$ is a $TK_5$ in $G'$ with branch vertices $x_1,x_2,y_1,y_2,z_2$.
If $w\in V(Q-a)$ then $z_2x_2\cup z_2Xy_2\cup (P_1\cup qQw\cup W)\cup (P_2\cup P\cup cCy_1)\cup (A\cup z_1Xy_2)\cup 
G[\{x_1, x_2, y_1,y_2\}]$ is a $TK_5$ in $G'$ with branch vertices $x_1, x_2, y_1, y_2, z_2$.
This completes the proof of (7).

\begin{itemize}
\item [(8)] We may assume that $G'$ has no path from $x_1Xz_1-x_1$ to $y_1$ and internally disjoint from $K\cup X$.
\end{itemize}
For, suppose that $R$ is a path in $G'$ from some $x\in V(x_1Xz_1-x_1)$ to $y_1$ and internally disjoint from $K\cup X$. 
Then $x\ne z_1$; as otherwise $z_2Bq\cup Q\cup aAz_1\cup R\cup y_1Cc\cup P\cup pBy_2$ is a path in $H$ through $z_2,z_1,y_1,y_2$ 
in order, contradicting (1). Let $P_1,P_2$ be the paths in (4) with $q^*=p$. 
We use $W$ from (7).  If $w\in V(A)-\{z_1, y_1\}$ then $z_2x_2\cup z_2Xy_2\cup 
(P_1\cup Q\cup aAw\cup W)\cup (P_2\cup P\cup cCy_1)\cup (R\cup xXy_2)\cup 
G[\{x_1,x_2,y_1,y_2\}]$ is a $TK_5$ in $G'$ with branch vertices $x_1,x_2,y_1,y_2,z_2$.
If $w\in V(C)-\{z_1, y_1\}$ then $z_2x_2\cup z_2Xy_2\cup (P_1\cup Q\cup aAy_1)\cup (P_2\cup P\cup cCw\cup W) \cup (R\cup xXy_2)\cup 
G[\{x_1,x_2,y_1,y_2\}]$ is a $TK_5$ in  $G'$ with branch vertices $x_1,x_2,y_1,y_2,z_2$.
This completes the proof of (8).

\begin{itemize}
\item [(9)] If $G'$ has a path from $x_1Xz_1-\{x_1, z_1\}$ to  $cCy_1-c$ and internally disjoint from $K\cup X$, 
then we may assume that  
\begin{itemize}
\item [$\bullet$] $w\in V(C)-\{y_1, z_1\}$ for any choice of $W$ in (7), and 
\item [$\bullet$] $G'$ has no path from $x_2$ to $C-\{y_1,z_1\}$ and internally disjoint from $K\cup X$.  
\end{itemize}
\end{itemize}
Let $S$ be a path in $G'$  from some $s\in V(x_1Xz_1)-\{x_1, z_1\}$ to  $V(cCy_1-c)$ and internally disjoint from $K\cup X$. 
Since $X$ is induced in $G'-x_1x_2$, $G'[H-\{y_2,z_1,z_2\}+s]$ is 2-connected.
Hence, since $N(x_1Xz_1-\{x_1,z_1\})\subseteq V(C')\cup \{x_1,z_1\}$ (by (6)),  
$G'$ has independent paths $S_1, S_2$ from $s$ to distinct $s_1, s_2\in V(C)-\{z_1, y_1\}$ 
and internally disjoint from $K\cup X$. Because of $S$, we may assume that $z_1, s_1, s_2, y_1$ occur on $C$ in this order and $s_2\in V(cCy_1-c)$.

Suppose we may choose the $W$  in (7) with $w\in V(A)-\{z_1, y_1\}$.  Let $P_1,P_2$ be the paths in (4) with $q^*=p$. 
Then $z_2x_2 \cup z_2Xy_2\cup sXx_1\cup sXy_2 \cup (P_2\cup P \cup cCs_1\cup S_1)\cup (S_2\cup s_2Cy_1\cup y_1x_2)\cup (P_1\cup Q\cup aAw\cup W) 
\cup G[\{x_1, x_2, y_2\}]$ is a $TK_5$ in $G'$ with branch vertices $s,x_1, x_2, y_2, z_2$. 

Now assume that $S'$ is a path in $G'$ from $x_2$ to some $s'\in V(C)-\{y_1,z_1\}$ and internally disjoint from $K\cup X$. 
Then $S_1\cup S_2\cup S'\cup (C-z_1)$ contains independent paths $S_1',S_2'$  which are from $s$ to $y_1,x_2$, respectively (when $s'\in V(z_1Cs_2)-\{s_2,z_1\}$),
 or from $s$ to $c,x_2$, respectively (when $s'\in V(s_2Cy_1-y_1)$). If $S'_1,S'_2$ end at $y_1,x_2$, respectively, 
then $sXx_1\cup sXy_2\cup S_1'\cup S_2'\cup (y_1Aa\cup Q\cup qBy_2)\cup G[\{x_1, x_2, y_1, y_2\}]$ 
is a $TK_5$ in $G'$ with branch vertices $s, x_1, x_2, y_1,y_2$.  
So assume that  $S_1'$, $S_2'$ end at $c, x_2$, respectively. Let $P_1,P_2$ be the paths in (4) with $q^*=p$.  
Then $sXx_1\cup sXy_2\cup z_2x_2\cup z_2Xy_2\cup (S_1'\cup P\cup P_2)\cup 
S_2'\cup  (P_1\cup Q\cup aAy_1 \cup y_1x_1)\cup G[\{x_1, x_2, y_2\}]$ is a $TK_5$ in 
 $G'$ with branch vertices $s,x_1, x_2, y_2, z_2$. This completes the proof of (9).

\medskip

The next two claims deal with $L(A)$ and $L(C)$. First, we may assume that 
\begin{itemize}
\item [(10)] $L(A)\cap A\subseteq z_1Aa$.
\end{itemize}
For any $(A\cup C)$-bridge $R$ of $H$ contained in $L(A)$, let $z(R), y(R)\in V(R\cap A)$ such that $z(R)Ay(R)$ is maximal.
Suppose for some $(A\cup C)$-bridge  $R_1$ of $H$ contained in $L(A)$, we have $y(R_1)Az(R_1)\not\subseteq z_1Aa$. 
Let $R_1, \ldots, R_m$ be a maximal sequence of $(A\cup C)$-bridges of $H$ contained in $L(A)$, 
such that  for each $i\in \{2,\ldots, m\}$,  $R_i$ contains an internal vertex of 
$\bigcup_{j=1}^{i-1}z(R_j)Ay(R_j)$ (which is a path). 
Let $a_1, a_2\in V(A)$ such that $\bigcup_{j=1}^{m}z(R_j)Ay(R_j)=a_1Aa_2$. 
By (c), $J(A,C)$ does not intersect $a_1Aa_2-\{a_1,a_2\}$; so  $a_1,a_2\in V(aAy_1)$. By (d), $G'$ has no path from $a_1Aa_2-\{a_1,a_2\}$ to $C$ and 
internally disjoint from $K\cup X$. Hence by (5), 
$\{a_1,a_2, x_1, x_2, y_2\}$ is a cut in $G$.  Thus,  $G$ has a separation $(G_1,G_2)$ such that 
$V(G_1\cap G_2)=\{a_1, a_2, x_1, x_2, y_2\}$, $P\cup Q\cup B'\cup C\cup X\subseteq G_1$, and $a_1Aa_2\cup \left(\bigcup_{j=1}^{m}R_j\right)
\subseteq G_2$.

Let $z\in V(G_2)-\{a_1,a_2, x_1, x_2, y_2\}$ and assume $z_1,a_1,a_2,y_1$ occur on $A$ in order. 
Since $G$ is 5-connected, $G_2-y_2$ contains four independent paths $R_1,R_2,R_3,R_4$ 
from $z$ to $x_1,x_2,a_1,a_2$, respectively. Now $R_1\cup R_2\cup (R_3\cup a_1Az_1\cup z_1Xy_2)\cup (R_4\cup a_2Ay_1)
\cup (y_1Cc\cup P\cup pBy_2)\cup G[\{x_1,x_2,y_1,y_2\}]$
is a $TK_5$ in $G'$ with branch vertices $x_1,x_2,y_1,y_2,z$.  
This completes the proof of (10).

\begin{itemize}
\item [(11)] We may assume that if $R$ is an $(A\cup C)$-bridge of $H$ contained in $L(C)$ and $R\cap (cCy_1-c)\ne \emptyset$ then $|V(R)-V(C)|=1$ and $N(R-C)=\{c_1,c_2,s_1,s_2,y_2\}$, 
with $c_1Cc_2=c_1c_2$ and $s_1s_2=s_1Xs_2\subseteq z_1Xx_1$.
\end{itemize}
For any $(A\cup C)$-bridge $R$ in $L(C)$, let $z(R),y(R)\in V(C\cap R)$ such that $z(R)Cy(R)$ is maximal.
Let $R_1$ be an  $(A\cup C)$-bridge of $H$ contained in $L(C)$ such that $R_1\cap (cCy_1-c)\ne \emptyset$. 

Let $R_1, \ldots, R_m$ be a maximal sequence of $(A\cup C)$-bridges of $H$ contained in $L(C)$, such that 
for each $i\in \{2,\ldots, m\}$,  $R_i$ contains an internal vertex of $\bigcup_{j=1}^{i-1}z(R_j)Cy(R_j)$ (which is a path). 
Let $c_1, c_2\in V(C)$ such that $c_1Cc_2=\bigcup_{j=1}^{m} z(R_j)Cy(R_j)$, with $z_1,c_1,c_2,y_1$ on $C$ in order. 
So   $c_2\in V(cCy_1-y_1)$ and, hence, $c_1\in V(cCy_1 - y_1)$ by (c) and the existence of $P$. Let $R'=   \bigcup_{j=1}^{m} R_j\cup c_1Cc_2$. 
 
By (c), $G'$ has no path from $c_1Cc_2-\{c_1,c_2\}$ to $V(B'\cup P\cup Q)\cup \{z_1\}$ and internally disjoint from $K\cup X$. By (d), $G'$ has no path from $c_1Cc_2-\{c_1,c_2\}$ 
to $A-\{y_1,z_1\}$ and internally disjoint from $K\cup X$.  

If $N(x_2)\cap V(R'-\{c_1,c_2\})\ne \emptyset$ then by (5) and (9), $N(R'-\{c_1,c_2\})=\{x_1,x_2,y_2,c_1,c_2\}$. 
Let $z\in V(R')-\{x_1,x_2,c_1,c_2\}$. 
Since $G$ is 5-connected, $R'$ has independent paths $W_1,W_2,W_3,W_4$ from $z$ to $x_1,x_2,c_2,c_1$, respectively. Now $W_1\cup W_2\cup (W_3 \cup c_2Cy_1)
\cup (W_4\cup c_1Cz_1\cup z_1Xy_2)\cup (y_1Aa\cup Q\cup qBy_2)\cup G[\{x_1,x_2,y_1,y_2\}]$ is a 
$TK_5$ in $G'$ with branch vertices $x_1,x_2,y_1,y_2,z$. 

So we may assume $N(x_2)\cap V(R'-\{c_1,c_2\})= \emptyset$. Since $G$ is 5-connected, it follows from (5) that  there exist distinct
$s_1, s_2\in V(x_1Xz_1-z_1)\cap N(R'-\{c_1, c_2\})$. Choose $s_1,s_2$  such that $s_1Xs_2$ is maximal and assume that $x_1, s_1, s_2, z_1$ 
occur on $X$ in this order.   By (6), $\{c_1, c_2, s_1, s_2,y_2\}$ is a $5$-cut in $G$; so $G$ has a separation $(G_1,G_2)$ 
such that $V(G_1\cap G_2)=\{c_1, c_2, s_1, s_2,y_2\}$ and $R'\cup c_1Cc_2\cup s_1Xs_2\subseteq G_2$. By (6) again, 
$(G_2-y_2, c_1,c_2,s_1,s_2)$ is planar (since $G$ is 5-connected). If $|V(G_2)|\ge 7$ then by Lemma~\ref{apexside1}, $(i)$ or $(ii)$ or $(iii)$ holds.
 So we may assume that $|V(G_2)|=6$, and we have the assertion of (11).  

\medskip
 We may assume that 
\begin{itemize}
\item [(12)] $H$ has  a path $Q'$ from $y_1$ to some $q'\in V(Q-a)$ and internally disjoint from $K$. 
\end{itemize}
First, suppose that $y_1\in V(J(A,C))$.
Then, $H$ has a path $Q'$ from $y_1$ to some $q'\in V(P-c)\cup V(Q-a)\cup V(B)$ internally disjoint from $K$.  We may assume $q'\in V(P-c)\cup V(B)$; 
for otherwise, $q'\in V(Q-a)$ and the claim holds. 
If $q'\in V(P-c)\cup V(y_2Bq-q)$ then $(P-c)\cup (y_2Bq-q)\cup Q'$ contains a path $Q''$ from $y_1$ to $y_2$; so 
$z_1Xx_1\cup z_1Xy_2\cup C\cup (z_1Aa\cup Q\cup qBz_2\cup z_2x_2)\cup Q''\cup G[\{x_1, x_2, y_1, y_2\}]$ is a $TK_5$ in $G'$ 
with branch vertices $x_1,x_2,y_1,y_2,z_1$.  Hence, we may assume $q'\in V(qBz_2-q)$. Let $P_1,P_2$ be the paths in (4) with $q^*=q'$. Then $z_2x_2\cup z_2Xy_2\cup 
(P_1\cup Q\cup aAz_1\cup z_1Xx_1)\cup (P_2\cup Q') \cup (y_1Cc\cup P\cup pBy_2)\cup G[\{x_1,x_2,y_1,y_2\}]$ is a $TK_5$ in $G'$ with branch vertices $x_1,x_2,
y_1,y_2,z_2$.

Thus, we may assume that $y_1\notin V(J(A,C))$. Note that $y_1\notin V(L(A))$ (by (10)) and $y_1\notin V(L(C))$ (by (8) and (11)).   
Hence, since $y_1y_2\notin E(G)$ and $G$ is 5-connected, $y_1$ is contained in some $(A\cup C)$-bridge of $H$, say $D_1$, 
with $D_1\subseteq L(A,C)$ and $D_1\ne J(A, C)$. Note that $|V(D_1)|\ge 3$ as $A$ and $C$ are induced paths. 
For any  $(A\cup C)$-bridge $D$ of $H$ with that  $D\subseteq L(A,C)$ and $D\ne J(A,C)$, 
let $a(D) \in V(A)\cap V(D)$ and $c(D)\in V(C)\cap V(D)$ such that $z_1Aa(D)$ and  $z_1Cc(D)$ are minimal. 

Let $D_1, \ldots, D_k$ be a maximal sequence of $(A\cup C)$-bridges of $H$ with $D_i\subseteq L(A,C)$ (so $D_i\ne J(A,C)$) for $i\in [k]$, 
such that, for each $i\in [k-1]$, $D_{i+1}\cap (A\cup C)$ is not contained in $\bigcup_{j=1}^i(c(D_j)Cy_1\cup a(D_j)Ay_1)$, and  
$D_{i+1}\cap (A\cup C)$ is not contained  in $\bigcap_{j=1}^i(z_1Cc(D_j) \cup z_1Aa(D_j))$. 
Note that for any $i\in [k]$, $\bigcup_{j=1}^ia(D_j)Ay_1$ and $\bigcup_{j=1}^ic(D_j)Cy_1$ are paths. So let $a_i\in V(A)$ and $c_i\in V(C)$ such that 
$\bigcup_{j=1}^ia(D_j)Ay_1=a_iAy_1$ and $\bigcup_{j=1}^ic(D_j)Cy_1=c_iCy_1$. Let $S_i=a_iCy_1\cup c_iCy_1\cup \left(\bigcup_{j=1}^iD_j\right)$.

Next, we claim that for any $l\in [k]$ and for any $r_l\in V(S_l)-\{a_l,c_l\}$ there exist three independent paths
$A_l,C_l,R_l$  in $S_l$ from $y_1$ to $a_l,c_l,r_l$, respectively.
This is clear  when $l=1$; note that if $a_l=y_1$, or $c_l=y_1$, or
$r_l=y_1$  then $A_l$, or $C_l$, or $R_l$  is a trivial path. 
Now assume that the assertion is  true for some $l\in [k-1]$. Let $r_{l+1}\in
V(S_{l+1})-\{a_{l+1},c_{l+1}\}$. When $r_{l+1}\in V(S_l)-\{a_l,c_l\}$ let
$r_l:=r_{l+1}$; otherwise, let $r_l\in V(D_{l+1})$ with $r_l\in V(a_lAy_1-a_l)\cup
V(c_lCy_1-c_l)$. By induction hypothesis, there
are three independent paths $A_l,C_l,R_l$ in $S_l$ from $y_1$ to
$a_l,c_l,r_l$, respectively. 
If $r_{l+1}\in V(S_l)-\{a_l,c_l\}$ then $A_{l+1}:=A_l\cup
a_lAa_{l+1}, C_{l+1}:=C_l\cup c_lCc_{l+1}, R_{l+1}:=R_l$ are the
desired paths in $S_{l+1}$. 
If $r_{l+1}\in V(D_{l+1})-V(A\cup C)$ then
let $P_{l+1}$ be a path in $D_{l+1}$ from $r_l$ to $r_{l+1}$ and 
internally disjoint from $A\cup C$; we see that $A_{l+1}:=A_l\cup
a_lAa_{l+1}, C_{l+1}:=C_l\cup c_lCc_{l+1}, R_{l+1}:=R_l\cup P_{l+1}$
are the desired paths in $S_{l+1}$. 
So we may assume by symmetry
that $r_{l+1}\in V(a_{l+1}Aa_l-a_{l+1})$. Let $Q_{l+1}$ be a path in
$D_{l+1}$ from $r_l$ to $a_{l+1}$ and  internally disjoint from $A\cup
C$. Now $R_{l+1}:=A_l\cup a_lAr_{l+1}, C_{l+1}:=C_l\cup c_lCc_{l+1},
A_{l+1}:=R_l\cup Q_{l+1}$ are the desired paths in $S_{l+1}$.

We claim that  $J(A, C)$ has no vertex in $(a_kAy_1\cup c_kCy_1) - \{a_k, c_k\}$.
For, suppose there exists $r\in V(J(A,C))$ such that $r\in V(a_kAy_1-a_k)\cup V(c_kCy_1-c_k)$. Then let $A_k,C_k,R_k$ be independent (induced) paths 
in $S_k$ from $y_1$ to $a_k,c_k,r$, respectively. Let $A',C'$ be obtained from $A,C$ by replacing $a_kAy_1,c_kCy_1$ with $A_k,C_k$, respectively. We see 
that $J(A',C')$ contains $J(A,C)$ and $r$, contradicting (c).

Therefore, $a\in V(z_1Aa_k)$ and $c\in V(z_1Cc_k)$. 
Moreover, no $(A\cup C)$-bridge of $H$ in $L(A)$ intersects  $a_kAy_1 - a_k$ (by (10)).
Let $S_k'$ be the union of $S_k$ and all $(A\cup C)$-bridges of $H$ contained in $L(C)$ and intersecting $c_kCy_1-c_k$.
Then by (5) and (11),  $N(S_k'-\{a_k, c_k\})-\{a_k, c_k,x_2,y_2\}\subseteq V(x_1Xz_1)$. Since  $G$ is $5$-connected, $N(S_k'-\{a_k, c_k\})-\{a_k, c_k,x_2,y_2\}\ne \emptyset$.

We may assume that $N(S_k'-\{a_k, c_k\})-\{y_2, x_2, a_k, c_k\}\ne \{x_1\}$. For, otherwise, 
$G$ has a separation $(G_1,G_2)$ such that $V(G_1\cap G_2)=\{a_k,c_k,x_1, x_2, y_2\}$ and $X\cup P\cup Q\subseteq G_1$, and $S_k'\subseteq G_2$. 
Clearly, $|V(G_1)|\ge 7$. Since $G$ is 5-connected and $y_1y_2\notin E(G)$,  $|V(G_2)|\ge 7$.  
Hence, the assertion follows from Lemma~\ref{5cut_triangle}.

Thus, we may let $z\in N(S_k'-\{a_k, c_k\})-\{a_k, c_k,x_1,x_2,y_2\}$ such that $x_1Xz$ is maximal.  
Then $z\ne z_1$. For otherwise, let $r\in V(S_k')-\{a_k, c_k\}$ such that $rz_1\in E(G)$. Let $r'=r$ if $r\in V(S_k)$ and, otherwise, 
let $r'\in V(c_kCy_1-c_k)$ with $r'r\in E(G)$ (which exists by (11)). 
Let $A_k,C_k,R_k$ be independent (induced) paths 
in $S_k$ from $y_1$ to $a_k,c_k,r'$, respectively. Now $z_2Bq\cup Q\cup aAz_1\cup (z_1rr'\cup R_k)\cup C_k\cup c_kCc\cup P\cup pBy_2$ 
is a path in $H$ through $z_2,z_1,y_1,y_2$ in order, contradicting (1). 

Let $C^*$ be the subgraph of $G$ induced by the union of $x_1Xz-x_1$ and the vertices of $L(C)-C$ adjacent to 
$c_kCy_1-c_k$ (each of which, by (11), has exactly two neighbors on $C$ and exactly two on $x_1Xz_1$).
Clearly, $C^*$ is connected.
 Let $G_z=G[x_1Xz\cup S_k'+x_2]$ and let $G_z'$ be the graph obtained from $G_z-\{x_1,x_2\}$ by contracting $C^*$ to a new vertex $c^*$. 

Note that $G_z'$ has no disjoint paths from $a_k, c_k$ to $c^*, y_1$, respectively; as otherwise, such paths, $c_kCc\cup P\cup pBy_2$, and $a_kAa\cup Q\cup qBz_2$
give two disjoint paths in $H$  which would contradict the choice of $Y,Z$. 
Hence, by Lemma~\ref{2path}, there exists a collection $\mathcal{A}$ of subsets of $V(G_z')- \{a_k, c_k, c^*, y_1\}$ such that 
$(G_z', \mathcal{A}, a_k, c_k, c^*, y_1)$ is $3$-planar. We choose $\mathcal{A}$ so that each member of $\mathcal{A}$ is minimal and, subject to this, 
$|{\cal A}|$ is minimal. 

We claim that  $\mathcal{A}=\emptyset$. For, let $T\in \mathcal{A}$. By (10), $T\cap V(L(A))=\emptyset$. Moreover,  $T\cap V(L(C))=\emptyset$; for 
otherwise, by (11), $c^*\in N(T)$ and $|N(T)\cap V(C)|=2$; 
so by (11) again (and since $C$ is induced in $H$),
$(G_z', \mathcal{A}-\{T\}, a_k, c_k, c^*, y_1)$ is 3-planar, contradicting the choice of  $\mathcal{A}$. 
Thus, $G[T]$ has a component, say $T'$, such that $T'\subseteq L(A,C)$. Hence, for any $t\in V(T')$, $L(A,C)$ has a path from $t$ to $aAy_1-y_1$  
(respectively, $cCy_1-y_1$) and internally disjoint from $A\cup C$. 
Since $G$ is 5-connected, $\{x_1,x_2\}\cap N(T')\ne \emptyset$. Therefore,  for some $i\in [2]$, $G'$ contains a path from $x_i$ to $aAy_1-y_1$ as well as a path from $x_i$ to 
$cCy_1-y_1$, both 
internally disjoint from $K\cup X$. However, this contradicts (9).

Hence,  $(G_z', a_k, c_k, c^*, y_1)$ is planar. So by (6) and (11), $(G_z-x_2, a_k, c_k, z, x_1,y_1)$ is planar.  
By (9) and (10), $N(x_2)\cap V(S_k)\subseteq V(a_kAy_1)$. Therefore, since $(G_z-x_2)-a_kAy_1$ is connected (by (10)),  
 $(G_z, a_k, c_k, z, x_2)$ is planar.

We claim that  $\{a_k, c_k, z, x_2, y_2\}$ is a $5$-cut in $G$. For, otherwise, by (7) and (9), 
$G'$ has a path $S_1$ from $x_1$ to $z_1Cc_k-\{z_1,c_k\}$ and internally 
disjoint from $K\cup X$. However,  $G'$ has a path $S_2$ from $z$ to $c_kXy_1-c_k$ and internally disjoint from $K\cup X$. 
Now $S_1,S_2$ contradict the second part of (6). 

Hence, $G$ has a separation $(G_1,G_2)$ such that $V(G_1\cap G_2)=\{a_k, c_k, z, x_2, y_2\}$, 
$B'\cup P\cup Q\cup X\subseteq G_1$, and $G_z\subseteq G_2$. Clearly, $|V(G_i)|\ge 7$ for $i\in [2]$. So 
 $(i)$ or $(ii)$ or $(iii)$ follows from Lemma~\ref{apexside1}.

\medskip

Now that we have established (12), the remainder of this proof will make heavy use of $Q'$. Our next goal is to obtain structure around $z_1$, which is done 
using claims (13) -- (17). We may assume  that
\begin{itemize}
\item [(13)] $x_1z_1\in E(X)$,  $w\in V(A) -\{y_1, z_1\}$ for any choice of $W$ in (7), and  $G'$ has no path from $x_2$ to $(A\cup C)-y_1$ and internally 
disjoint from $K\cup Q'\cup X$.
\end{itemize}
Let $P_1$, $P_2$ be the paths in (4) with $q^*=p$. 
Suppose $x_1z_1\notin E(X)$. Let $x_1s\in E(X)$. By (6), $G$ has a path $S$ from $s$ to some $s'\in V(C)-\{y_1,z_1\}$ and 
internally disjoint from $K\cup Q'\cup X$ (as $Q'\subseteq J(A,C)$). Hence,  
$z_2x_2\cup z_2Xy_2\cup (P_1\cup qQq'\cup Q')\cup ( P_2\cup P\cup cCs'\cup S\cup sx_1)\cup (A\cup z_1Xy_2)\cup G[\{x_1, x_2, y_1, y_2\}]$ 
is a $TK_5$ in  $G'$ with branch vertices $x_1,x_2,y_1,y_2,z_2$. 

Now suppose $W$ is a path in (7) ending at $w\in V(C)-\{y_1,z_1\}$. Then 
$z_2x_2\cup z_2Xy_2\cup  (P_1\cup qQq'\cup Q')\cup (P_2\cup P\cup cCw\cup W)\cup (A\cup z_1Xy_2)\cup G[\{x_1, x_2, y_1,y_2\}]$ is a $TK_5$ in 
$G'$ with branch vertices $x_1,x_2,y_1,y_2,z_2$.

Finally, suppose $G'$ has a path  $S$ from $x_2$ to some  $s\in V(A\cup C)-\{y_1\}$ and internally disjoint from $K\cup Q'\cup X$.
If $s\in V(A-y_1)$ then $z_1x_1\cup z_1Xy_2\cup C\cup (z_1As\cup S)\cup (Q'\cup q'Qq\cup qBy_2)\cup G[\{x_1, x_2, y_1, y_2\}]$ is a 
$TK_5$ in   $G'$ with branch vertices $x_1,x_2,y_1,y_2,z_1$. 
If $s\in V(C-y_1)$ then $z_1x_1\cup z_1Xy_2\cup A\cup (z_1Cs\cup S)\cup (Q'\cup q'Qq\cup qBy_2)\cup G[\{x_1, x_2, y_1, y_2\}]$ is a $TK_5$ in 
$G'$ with branch vertices $x_1,x_2,y_1,y_2,z_1$.

\begin{itemize}
\item [(14)] We may assume that $G'$ has no path from  $y_2Xz_2$ to  $(A\cup C)-y_1$ and  internally disjoint from $K\cup Q'\cup X$, 
and no path from $y_2Xz_1-z_1$ to $A-z_1$ and internally disjoint from $K\cup Q'\cup X$.
\end{itemize}
First, suppose $S$ is a path in $G'$  from  some $s\in V(y_2Xz_2)$ to some $s'\in V(A\cup C)-\{y_1\}$ and internally disjoint from $K\cup Q'\cup X$. 
Then $s\ne y_2$ as $N_{G'}(y_2)=\{w_1,w_2,w_3,x_1,x_2\}$. 
If $s'\in V(C-y_1)$ then $z_1x_1\cup z_1Xy_2\cup A\cup (z_1Cs'\cup S\cup sXx_2)\cup (Q'\cup q'Qq\cup qBy_2)\cup G[\{x_1, x_2, y_1, y_2\}]$ is a 
$TK_5$ in $G'$ with branch vertices $x_1,x_2,y_1,y_2,z_1$. 
If $s'\in V(A-y_1)$ then $z_1x_1\cup z_1Xy_2\cup C\cup (z_1As'\cup S\cup sXx_2)\cup (Q'\cup q'Qq\cup qBy_2)\cup G[\{x_1, x_2, y_1, y_2\}]$ is a $TK_5$ 
 in $G'$ with branch vertices $x_1,x_2,y_1,y_2,z_1$.

Now suppose  $S$ is a path in $G'$ from $s\in V(y_2Xz_1-z_1)$ to $s'\in V(A-z_1)$ and internally disjoint from $K\cup Q'\cup X$. 
Let $P_1,P_2$ be the paths in  (4) with $q^*=p$. Then $z_2x_2\cup z_2Xy_2\cup (P_1\cup qQq'\cup Q')\cup 
(P_2\cup P\cup cCz_1\cup z_1x_1)\cup (y_1As'\cup S\cup sXy_2) \cup G[\{x_1, x_2, y_1, y_2\}]$ is a $TK_5$  in $G'$ with branch vertices $x_1,x_2,y_1,y_2,z_2$.

\begin{itemize}
\item [(15)]  We may assume that 
\begin{itemize}
\item [$\bullet$] $J(A, C)\cap (z_1Cc-c)=\emptyset$, 
\item [$\bullet$] any path in $J(A, C)$ from $A-\{y_1, z_1\}$ to $(P-c)\cup (Q-a)\cup (Q'-y_1)\cup B$ and internally disjoint from $K\cup Q'$  must end on $(Q\cup Q')-q$, and 
\item  [$\bullet$] for any $(A\cup C)$-bridge $D$ of $H$ with $D\ne J(A,C)$, if $V(D)\cap V(z_1Cc-c)\ne \emptyset$ and $u\in V(D)\cap V(z_1Ay_1-z_1)$ 
 then $J(A, C)\cap (z_1Au-\{z_1, u\})=\emptyset$. 
\end{itemize}
\end{itemize}
First, suppose there exists $s\in V(J(A, C))\cap V(z_1Cc-c)$. Then $H$ has a path $S$ from $s$ to some $s'\in V(P-c)\cup V(Q-a)\cup V(Q'-y_1)\cup V(B-y_2)$ 
and internally disjoint from $K\cup Q'$.  If $s'\in V(Q'-y_1)\cup V(Q-a)\cup V(z_2Bp-p)$ then 
$S\cup (Q'-y_1)\cup (Q-a)\cup (z_2Bp-p)$ contains a path $S'$ from $s$ to $z_2$; 
so  $S'\cup sCz_1\cup A\cup y_1Cc\cup P\cup pBy_2$ is a path in $H$ through $z_2, z_1, y_1, y_2$ in order, contradicting (1).
Hence, $s'\in V(P-c)\cup V(y_2Bp-y_2)$ and, by (2), $s = z_1$. 
Let $P_1,P_2$ be the paths in (4)  with $q^*=p$ (if $s'\in V(P-c)$) or $q^*=s'$ (if $s'\in V(y_2Bp)-\{p,y_2\}$). 
Then $S\cup (P-c)\cup P_2$ contains a path $S'$ from $z_1$ to $z_2$. Let $W,w$ 
be given as in (7). By (13), $w\in V(A)-\{y_1,z_1\}$.  Now $z_2x_2\cup z_2Xy_2\cup z_1x_1\cup z_1Xy_2\cup S'\cup (P_1\cup Q\cup aAw\cup W)\cup (C\cup y_1x_2)\cup G[\{x_1, x_2, y_2\}]$ is a 
$TK_5$ in  $G'$ with branch vertices $x_1,x_2,y_2,z_1,z_2$.
 
Now suppose $S$ is path in  $J(A, C)$ from $s\in V(A-\{y_1, z_1\})$ to $s'\in V(P-c)\cup V(B-q)$ and internally disjoint from $K\cup Q'$. 
Since $N_{G'}(y_2)=\{w_1,w_2,w_3,x_1,x_2\}$, $s'\ne y_2$. Let $P_1,P_2$ be the paths 
in (4) with $q^*=p$ (if  $s'\in V(P-c)$) or $q^*=s'$ (if $s'\in V(B-q)$). Let $S'$ be a path in $P_2\cup S\cup (P-c)$ from $s$ to $z_2$. Let $W,w$ 
be given as in (7). By (13), $w\in V(A)-\{y_1,z_1\}$. Hence, 
$z_2x_2\cup z_2Xy_2\cup (P_1\cup qQq'\cup Q')\cup (S'\cup sAw\cup W)\cup (C\cup z_1Xy_2)\cup G[\{x_1, x_2, y_1, y_2\}]$ is a $TK_5$  in  $G'$ with branch  vertices $x_1,x_2,y_1,y_2,z_2$.

Finally, suppose $D$ is some $(A\cup C)$-bridge of $H$ with $D\ne J(A,C)$, $v\in V(D)\cap V(z_1Cc-c)$, and $u\in V(D)\cap V(z_1Ay_1-z_1)$. Then $D$ has a path $T$ 
from $v$ to $u$ and internally disjoint from $K\cup Q'$. If there exists $s\in V(J(A, C))\cap V(z_1Au-\{z_1, u\})$ then 
$J(A,C)$ has a path $S$ from $s$ to some $s'\in V(Q-a)$ and internally disjoint from $K$. 
Now $z_2Bq\cup qQs'\cup S\cup sAz_1\cup z_1Cv\cup T\cup uAy_1\cup y_1Cc\cup P\cup pBy_2$ is a path in $H$ through $z_2, z_1, y_1, y_2$ in order, contradicting  (1).

\begin{itemize}
\item [(16)] We may assume $L(A)=\emptyset$.
\end{itemize}
Suppose $L(A)\ne \emptyset$. For each $(A\cup C)$-bridge $R$ of $H$ contained in $L(A)$, let $a_1(R),a_2(R)\in V(R\cap A)$ with $a_1(R)Aa_2(R)$ 
maximal. Let $R_1, \ldots, R_m$ be a maximal sequence of $(A\cup C)$-bridges of $H$ contained in $L(A)$, such that for $i=2, \ldots, m$, 
  $R_i$ contains an internal vertex of  $\bigcup_{j=1}^{i-1}(a_1(R_j)Aa_2(R_j))$ (which is a path).
Let  $a_1, a_2\in V(A)$ such that  $\bigcup_{j=1}^{m}a_1(R_j)Aa_2(R_j)=a_1Aa_2$. Let $L=\bigcup_{j=1}^{m}R_j$.

 By (c), $J(A,C)\cap (a_1Aa_2-\{a_1,a_2\})=\emptyset$. 
By (d), $L(A, C)\cap (a_1Aa_2-\{a_1,a_2\})=\emptyset$. By (10), $a_1,a_2\in V(z_1Aa)$. So $z_1\notin N(L\cup a_1Aa_2-\{a_1,a_2\})$. 
Hence by (14),  $V(z_1Xz_2-y_2)\cap N(L\cup a_1Aa_2-\{a_1,a_2\})=\emptyset$.
By (13), $x_2\notin N(L\cup a_1Aa_2-\{a_1,a_2\})$.
Thus,  $\{a_1, a_2,x_1,y_2\}$ is a cut in $G$ separating $L$ from $X$, which is a contradiction (since $G$ is $5$-connected). 

\begin{itemize}
\item [(17)] $z_1c\in E(C)$, $z_1y_2\in E(G)$, and $z_1$ has degree 5 in $G$. 
\end{itemize}
Let $C^*$ be the union of $z_1Cc$ and  all $(A\cup C)$-bridges of $H$ intersecting $z_1Cc-c$.  By (15), $V(C^*\cap J(A,C))=\{c\}$. 

Suppose (17) fails.  If $C^*=z_1Cc$ then, since $A, C$ are induced paths and $L(A)=\emptyset$ (by (16)), $z_1y_2\in E(G)$ and $z_1Cc\ne z_1c$; so any vertex of $z_1Cc-\{c,z_1\}$ 
would have degree 2 in $G$ (by (15)), a contradiction. So $C^*-z_1Cc\ne \emptyset$. 
Since $G'-X$ is 2-connected,  $(C^*-z_1Cc)\cap (A-z_1)\ne \emptyset$ by (c) (and since $J(A.C)\cap \cap (zCc-c)=\emptyset$ by (15)). 
Moreover, if $|V(z_1Cc)|\ge 3$ then there is a path in $C^*$ from $z_1Cc-\{c,z_1\}$ to $A-z_1$ and internally disjoint from $A\cup C$. 

Let $a^*\in V(A\cap C^*)$ with  $a^*Ay_1$ minimal, and let $u\in V(z_1Xy_2)$ with $uXy_2$ minimal such that $u$ 
is a neighbor of $(C^*-c)\cup (z_1Aa^*-a^*)$. 

We may assume that $\{a^*,c,u,x_1, y_2\}$ is a $5$-cut in $G$. First, note, by (15), that $J(A, C)\cap ((z_1Aa^*-a^*)\cup (z_1Cc-c)) =\emptyset$ (in particular, 
$a^*\in V(z_1Aa)$). 
Hence, if $u=z_1$ then it is clear from (d), (13) and (14) that  $\{a^*,c,u,x_1, y_2\}$ is a $5$-cut in $G$. So we may assume $u\ne z_1$. 
Then $G'$ contains a path $T$ from $u$ to $u'\in V(A-z_1)$ 
and internally disjoint from $A\cup cCy_1\cup P\cup Q\cup Q'\cup B'$. 
Suppose  $\{ a^*, c, u,x_1, y_2\}$ is not a $5$-cut in $G$. Then  by (d), (13) and (14), $G'$ has a path $R$ 
from $r\in V(z_1Xu-u)$ to $r'\in V(P-c)\cup V(Q-a)\cup V(Q'-y_1)\cup V(B')$ 
and internally disjoint from $K\cup X$. Note that $r'\ne y_2$ as $N_{G'}(y_2)=\{w_1,w_2,w_3,x_1,x_2\}$. 
If $r'\in V(B'-q)$ then let $P_1,P_2$ be the paths in (4) with $q^*=r'$; now $z_2x_2\cup z_2Xy_2\cup (P_1\cup qQq'\cup Q')\cup (P_2\cup R\cup rXx_1)\cup (y_1Au'\cup T\cup uXy_2)
\cup G[\{x_1,x_2,y_1,y_2\}]$ is a $TK_5$ in $G$ with branch vertices $x_1,x_2,y_1,y_2,z_2$. If $r'\in V(P-c)$ then   
let $P_1,P_2$ be the  paths in (4) with $q^*=p$; now $z_2x_2\cup z_2Xy_2\cup (P_1\cup qQq'\cup Q')\cup (P_2\cup pPr'\cup R\cup rXx_1)\cup (y_1Au'\cup T\cup uXy_2)
\cup G[\{x_1,x_2,y_1,y_2\}]$ is a $TK_5$ in $G$ with branch vertices $x_1,x_2,y_1,y_2,z_2$. Now assume $r'\in V(Q-a)\cup V(Q'-y_1)$. 
Then $(Q-a)\cup (Q'-y_1)\cup R$ contains a path $R'$ from $r$ to $q$.
Let  $P_1,P_2$ be the paths in (4) with $q^*=p$; now $z_2x_2\cup z_2Xy_2\cup (P_1\cup R'\cup rXx_1)\cup (P_2\cup P\cup cCy_1)\cup (y_1Au'\cup T\cup uXy_2)
\cup G[\{x_1,x_2,y_1,y_2\}]$ is a $TK_5$ in $G$ with branch vertices $x_1,x_2,y_1,y_2,z_2$. 

Thus, $G$ has a separation $(G_1,G_2)$ such that  $V(G_1\cap G_2)=\{a^*,c,u,x_1, y_2\}$, $uXx_2\cup P\cup Q\subseteq G_1$, and $C^*\cup z_1Cc\cup z_1Aa^*\subseteq G_2$. 
Suppose  $G_2-y_2$ contains disjoint paths $T_1,T_2$ from $u, x_1$ to $a^*, c$, respectively. 
Let $P_1,P_2$ be the paths in (4) with $q^*=p$. Then $z_2x_2\cup z_2Xy_2\cup (P_1\cup qQq'\cup Q')\cup (P_2\cup P\cup T_2)\cup (y_1Aa^*\cup T_1\cup uXy_2)
\cup G[\{x_1, x_2, y_1, y_2\}]$ is a $TK_5$ in  $G'$ with  branch vertices $x_1,x_2,y_1,y_2,z_2$.  
So we may assume that such $T_1,T_2$ do not exist. Then by Lemma~\ref{2path}, $(G_2-y_2,u,x_1,a^*,c)$ is planar (as $G$ is 5-connected). 
If $|V(G_2)|\ge 7$ then, by Lemma~\ref{apexside1}, $(i)$ or $(ii)$ or $(iii)$ holds. Hence, we may assume that $|V(G_2)|=6$ and, hence, we have (17).

\medskip

We have now forced a structure around $z_1$. Next, we study the structure of $G'[B'\cup y_2Xz_2]$ to complete the proof of Theorem~\ref{y_2}.
We may assume that 

\begin{itemize}
\item [(18)] $(G'[B'\cup y_2Xz_2], p, q, z_2, y_2)$ is $3$-planar. 
\end{itemize}
For, otherwise,  by Lemma~\ref{2path}, $G'[B'\cup y_2Xz_2]$ has disjoint paths $R_1, R_2$  from $q, p$ to $y_2, z_2$, respectively. 
Now $z_1x_1\cup z_1Xy_2\cup A\cup  (z_1Cc\cup P\cup R_2\cup z_2x_2)\cup (R_1\cup qQq'\cup Q')\cup G[\{x_1, x_2, y_1, y_2\}]$ is a 
$TK_5$ in  $G'$ with branch vertices $x_1,x_2,y_1,y_2,z_1$. So we may assume (18). 

\medskip

Since $G$ is 5-connected, $G$  is $(5, V(K\cup Q'\cup y_2Xx_2 \cup z_1x_1))$-connected. Recall that $w_1y_2\in E(x_1Xy_2)$. 
Then  $w_1y_2$ and $w_1Xz_1$ are independent paths in $G$ from $w_1$ to $y_2,z_1$, respectively. 
So by Lemma~\ref{perfect},  $G$ has  five independent paths $Z_1, Z_2, Z_3, Z_4, Z_5$ from $w_1$ to $z_1$, $y_2$, $z_3, z_4, z_5$,  respectively, and  
internally disjoint from $K\cup Q'\cup y_2Xx_2\cup z_1x_1$, where $z_3, z_4, z_5\in 
V(K\cup Q'\cup y_2Xx_2 \cup z_1x_1)$.
 Note that we may assume $Z_2=w_1y_2$. Hence, $Z_1,Z_2,Z_3,Z_4,Z_5$ are paths in $G'$. 
By the fact that $X$ is induced, by (14), and by (5) and (17),
$z_3, z_4, z_5\in V(P)\cup V(Q-a)\cup V(Q')\cup V(B'-y_2)$. Recall that $L(A)=\emptyset$ from (16), and recall $W$ and $w$ from (7) and (13).

\begin{itemize}
\item [(19)] We may assume that  at least two of  $Z_3, Z_4, Z_5$ end in $B'-y_2$.
\end{itemize}
First, suppose at least two of $Z_3, Z_4, Z_5$ end on $P$. Without loss of generality, let $c, z_3, z_4, p$ occur on $P$ in this order. Let $P_1,P_2$ be 
the paths in (4) with 
$q^*=p$. Then $(Z_1\cup z_1x_1)\cup Z_2\cup z_2x_2\cup z_2Xy_2\cup (Z_4\cup z_4Pp\cup P_2)\cup (Z_3\cup z_3Pc\cup cCy_1\cup y_1x_2)\cup 
(P_1\cup Q\cup aAw \cup W)\cup G[\{x_1, x_2, y_2\}]$ is a $TK_5$  in  $G'$ with branch vertices $w_1,x_1,x_2,y_2,z_2$.

Now assume at least two of $Z_3, Z_4, Z_5$ are on $Q\cup Q'$, say $Z_3$ and $Z_4$. Then $Z_3\cup Z_4\cup Q\cup Q'$ contains 
two independent paths $Z_3',Z_4'$ from $w_1$ to $z',q$, respectively, 
where $z'\in \{a,y_1\}$. Hence 
$(Z_1\cup z_1x_1)\cup Z_2\cup (Z_3'\cup z'Ay_1)\cup (Z_4'\cup qBz_2\cup z_2x_2)\cup (y_2Bp\cup P\cup cCy_1)\cup G[\{x_1, x_2, y_1, y_2\}]$ is a 
$TK_5$  in  $G'$ with branch vertices $w_1,x_1,x_2,y_1,y_2$. 

So we may assume that $z_3\in V(B')-\{p,q\}$, and hence $Z_3=w_1z_3$. 
Suppose none of $Z_4,Z_5$ ends in $B'-y_2$. Then we may assume  $z_4\in V(P-p)$. Let $P_1,P_2$ be the paths in (4) with $q^*=z_3$. Then 
$(Z_1\cup z_1x_1)\cup Z_2 \cup z_2x_2\cup z_2Xy_2\cup (Z_3\cup P_2)
\cup (P_1\cup Q\cup aAw\cup W)\cup (Z_4\cup z_4Pc\cup cCy_1\cup y_1x_2)\cup G[\{x_1,x_2,y_2\}]$ is a $TK_5$ in $G'$ with branch vertices $w_1, x_1,x_2,y_2, z_2$.

\begin{itemize}
\item [(20)] We may assume that 
\begin{itemize}
\item [$\bullet$] $w_1$ has at most one neighbor in $B'$ that is in $qBz_2$ or separated from $y_2Bp$ in $G'[B'\cup y_2Xz_2]$ by a 2-cut contained in 
$qBz_2$, and 
\item [$\bullet$] $w_1$ has at most one neighbor in $B'$ that is in $y_2Bp-y_2$ or separated from $qBz_2$ in $G'[B'\cup y_2Xz_2]$ by a 2-cut contained in $y_2Bp$.
\end{itemize}
\end{itemize}
Suppose  there exist distinct $v_1,v_2\in  N(w_1)\cap V(B')$ such that for $i\in [2]$, $v_i\in V(qBz_2)$ or 
$G'[B'\cup y_2Xz_2]$ has a 2-cut contained in $qBz_2$ and separating $v_i$ from $y_2Bp$. Then, since $(G'[B'\cup y_2Xz_2], p, q, z_2, y_2)$ is $3$-planar 
(by (18)) and $H-y_2$ is 2-connected, $G'[B'+w_1]-y_2Bp$ contains independent paths 
$S_1,S_2$ from $w_1$ to $q,z_2$, respectively. Now $w_1Xx_1\cup w_1y_2\cup (S_1\cup qQq'\cup Q')\cup (S_2\cup z_2x_2)\cup (y_1Cc\cup P\cup pBy_2)\cup 
G[\{x_1,x_2,y_1,y_2\}]$ is a $TK_5$ in $G'$ with branch vertices $w_1,x_1,x_2,y_1,y_2$. 
 
Now suppose there exist distinct $v_1,v_2\in  N(w_1)\cap V(B')$ such that for $i\in [2]$, $v_i\in V(y_2Bp)$ or 
$G'[B'\cup y_2Xz_2]$ has a 2-cut contained in $y_2Bp$ and separating $v_i$ from $qBz_2$.
Then, since $(G'[B'\cup y_2Xz_2], p, q, z_2, y_2)$ is $3$-planar (by (18)) and $H-y_2$ is 2-connected, 
  $G'[B'+w_1]-(qBz_2-z_2)$ has independent paths $S_1,S_2$ from $w_1$ to $p,z_2$, respectively. Now 
$w_1Xx_1\cup w_1y_2\cup z_2x_2\cup z_2Xy_2\cup S_2\cup (S_1\cup P\cup cCy_1\cup y_1x_2)\cup (z_2Bq\cup Q\cup aAw\cup W) \cup G[\{x_1,x_2,y_2\}]$ is a 
$TK_5$ in $G'$ with branch vertices $w_1,x_1,x_2,y_2,z_2$.

\begin{itemize}
\item [(21)] $G'[B'\cup y_2Xz_2]$ has a 2-separation $(B_1,B_2)$ such that  $N(w_1)\cap V(B'-y_2)\subseteq V(B_1)$, $pBq\subseteq B_1$, and $y_2Xz_2\subseteq B_2$.
\end{itemize}
Let $z \in N(w_1)\cap V(B')$ be arbitrary. If there exists a path $S$ in $B'-(pBy_2\cup (qBz_2-z_2))$ from $z_2$ to $z$ then 
$z_2x_2\cup z_2Xy_2\cup (z_2Bq\cup qQq'\cup Q')\cup (S\cup zw_1\cup w_1Xx_1)\cup (y_1Cc\cup P\cup pBy_2)\cup  G[\{x_1, x_2, y_1, y_2\}]$ 
is a $TK_5$ in  $G'$ with branch vertices $x_1,x_2,y_1,y_2,z_2$. 
So we may assume that such path $S$ does not exist. Then, since $(G'[B'\cup y_2Xz_2], p,q,z_2,y_2)$ is 3-planar (by (18)) and $G'-X$ is 2-connected, 
$z\in V(y_2Xp\cup qBz_2)$ (in which case let $B_z'=z$ and $B_z''=G'[B'\cup y_2Xz_2]$), or   
$G'[B'\cup y_2Xz_2]$ has a 2-separation  $(B_z',B_z'')$ such that $B_z'\cap B_z''\subseteq y_2Bp\cup qBz_2\cup y_2Xz_2$, $z\in V(B_z'-B_z'')$ and $z_2\in V(B_z''-B_z')$.

We claim that we may assume that $w_1$ has exactly two neighbors in $B'$, say $v_1,v_2$, such that $v_1\in V(y_2Bp-y_2)$ or $G'[B'\cup y_2Xz_2]$ has a 2-cut contained in $y_2Bp$ 
and separating $v_1$ from $qBz_2$, and  $v_2\in V(qBz_2-z_2)$ or $G'[B'\cup y_2Xz_2]$ has a 2-cut contained in $qBz_2$ and separating $v_2$ 
from $y_2Bp$. This follows from (20) if for every choice of $z$, $B_z'\cap B_z''\subseteq y_2Bp$ or  $B_z'\cap B_z''\subseteq qBz_2$.
So we may assume that there exists $v\in  N(w_1)\cap V(B')$ such that $pBq\subseteq B_v'$ and  we choose $v$ and $(B_v',B_v'')$ with $B_v'$ maximal. 
If $pBq\subseteq B_z'$ for all choices of $z$ then, by (18), we have (21). Thus, we may 
assume that there exists $z\in  N(w_1)\cap V(B')$ such that $pBq\not\subseteq B_z'$ for any choice of $(B_z',B_z'')$.  
Then $B_z'\cap B_z''\subseteq y_2Bp$ or  $B_z'\cap B_z''\subseteq qBz_2$. First, assume $B_z'\cap B_z''\subseteq qBz_2$. Then by the maximality of $B_v'$, 
$B'-y_2Bp$ has independent paths $T_1,T_2$ from $z_2$ to $q,z$, respectively. Hence, $z_2x_2\cup z_2Xy_2\cup (T_1\cup qQq'\cup Q')\cup (T_2\cup zw_1\cup w_1Xx_1)\cup (y_1Cc\cup P\cup pBy_2)
\cup G[\{x_1,x_2,y_1,y_2\}]$ is a $TK_5$ in $G'$ with branch vertices $x_1,x_2,y_1,y_2,z_2$.
Now assume $B_z'\cap B_z''\subseteq y_2Bp$. Then by (20), for any $t\in N(w_1)\cap V(B_v')$, $t\notin 
V(y_2Bp-y_2)$ and $G'[B'\cup y_2Xz_2]$ has no 2-cut contained in $y_2Bp$ 
and separating $t$ from $qBz_2$. If for every choice of $t\in N(w_1)\cap V(B_v')$, 
we have $t\in V(qBz_2-z_2)$ or $G'[B'\cup y_2Xz_2]$ has a 2-cut contained in $qBz_2$ and separating $t$ 
from $y_2Bp$ then the claim follows from (20). Hence, we may assume that  $t$ can be chosen   so that $t\notin V(qBz_2-z_2)$ and 
$G'[B'\cup y_2Xz_2]$ has no 2-cut contained in $qBz_2$ and separating $t$ from $y_2Bp$. Then, by (18) and 2-connectedness of $G'-X$, $G[B'+w_1]-(qBz_2-z_2)$ has independent paths $S_1,S_2$ 
from $w_1$ to $p,z_2$, respectively. Now $w_1Xx_1\cup w_1y_2\cup z_2x_2\cup z_2Xy_2\cup S_2\cup (S_1\cup P\cup cCy_1\cup y_1x_2)\cup (z_2Bq\cup Q\cup aAw\cup W) \cup G[\{x_1,x_2,y_2\}]$ is a 
$TK_5$ in $G'$ with branch vertices $w_1,x_1,x_2,y_2,z_2$.

 Thus, we may assume that $Z_3=w_1v_1$, $Z_4=w_1v_2$, and $Z_5$ ends at some $v_3\in V(P\cup Q\cup Q')-\{a,p,q\}$. 
Suppose $v_3\in V(P-p)$. Let $P_1,P_2$ be the paths in (4) with $q^*=v_1$. Then $w_1Xx_1\cup w_1y_2\cup z_2x_2\cup z_2Xy_2
\cup (w_1v_1\cup P_2)\cup (Z_5\cup v_3Pc\cup cCy_1\cup y_1x_2)\cup (P_1\cup Q\cup aAw\cup W)\cup   G[\{x_1,x_2,y_2\}]$ is a 
$TK_5$ in $G'$ with branch vertices $w_1,x_1,x_2,y_2,z_2$. 

Now assume $v_3\in V(Q\cup Q')-\{a,q\}$. Then $(B'-y_2Bp)\cup Z_5\cup Q\cup Q'\cup (A-z_1)\cup w_1v_2$ has independent paths $R_1,R_2$ from $w_1$ to $y_1,z_2$, respectively. 
So $w_1Xx_1\cup w_1y_2\cup R_1\cup (R_2\cup z_2x_2)\cup (y_1Cc\cup P\cup pBy_2)\cup  G[\{x_1,x_2,y_1,y_2\}]$ is a 
$TK_5$ in $G'$ with branch vertices $w_1,x_1,x_2,y_1,y_2$.  This completes the proof of (21). 

\medskip

By (21), let $V(B_1\cap B_2)=\{t_1, t_2\}$ with $t_1\in V(y_2Bp)$ and $t_2\in V(qBz_2)$.  Choose $\{t_1,t_2\}$ so that $B_2$ is minimal.
Then we may assume that $(G'[B_2+x_2],t_1,t_2,x_2,y_2)$ is 3-planar. For, otherwise, by Lemma~\ref{2path}, $G'[B_2+x_2]$ contains disjoint paths $T_1,T_2$ from $t_1,t_2$ to 
$x_2,y_2$, respectively. Then $z_1x_1\cup z_1Xy_2\cup A\cup (z_1Cc\cup P\cup pBt_1\cup T_1)\cup (Q'\cup q'Qq\cup qBt_2\cup T_2) \cup 
  G[\{x_1, x_2, y_1, y_2\}]$ is a $TK_5$ in  $G'$ with branch vertices $x_1,x_2,y_1,y_2,z_1$.

Suppose there exists $ss'\in E(G)$ such that $s\in V(z_1Xw_1-w_1)$ and $s'\in  V(B_2)-\{t_1,t_2\}$. Then $s'\notin V(X)$, as $X$ is induced in $G'-x_1x_2$. 
By (19), (20) and (21), we may assume that $B_1-qBt_2$ contains a path $R$ from $z_3$ to $p$. 
By the minimality of $B_2$ and 2-connectedness of $H-y_2$, $(B_2-t_1)-(y_2Xz_2-z_2)$ contains independent paths $R_1,R_2$ from $z_2$ to $s',t_2$, respectively. 
Now $z_2x_2\cup z_2Xy_2\cup (R_1\cup s's\cup sXx_1)\cup (R_2\cup t_2Bq\cup qQq'\cup Q')\cup (y_1Cc\cup P\cup R\cup z_3w_1y_2)
\cup  G[\{x_1, x_2, y_1, y_2\}]$ is a $TK_5$ in  $G'$ with branch vertices $x_1,x_2,y_1,y_2,z_2$.

Thus, we may assume that $ss'$ does not exist. 
Since $G$ is 5-connected, $\{t_1,t_2,y_2,x_2\}$ is not a cut. So $H$ has a path $T$ from some $t\in V(y_2Xx_2)-\{y_2, x_2\}$ to some 
$t'\in V(P\cup Q\cup Q'\cup A\cup C)-\{p,q\}$ and internally disjoint 
from $K\cup Q'$. 
By (14),  $t'\notin V(A\cup C)-\{y_1\}$. 

If $t'\in V(P-p)$ then $z_1x_1\cup z_1Xy_2\cup A\cup (z_1Cc\cup cPt'\cup T \cup tXx_2)\cup 
(Q'\cup q'Qq\cup qBy_2)\cup  G[\{x_1, x_2, y_1, y_2\}]$ is a $TK_5$ in  $G'$ with branch vertices $x_1,x_2,y_1,y_2,z_1$.  
So we assume $t'\in V(Q\cup Q')-\{a,q\}$.
 
If $q\ne q'$ or $t'\in V(Q')$ then  $(T\cup Q\cup Q')-q$ has a path $Q^*$ from $t$ to $y_1$; now 
$z_1x_1\cup z_1Xy_2\cup A\cup (z_1Cc\cup P\cup pBz_2\cup z_2x_2)\cup (Q^*\cup  sXy_2)\cup  G[\{x_1, x_2, y_1, y_2\}]$ is 
a $TK_5$ in  $G'$ with branch vertices $x_1,x_2,y_1,y_2,z_1$. 
So assume $q=q'$ and $t'\in V(Q)-\{a,q\}$. 
Then $z_1x_1\cup z_1Xy_2\cup C\cup (z_1Aa\cup aQt'\cup T\cup tXx_2)
\cup (Q'\cup qBy_2)\cup G[\{x_1, x_2, y_1, y_2\}]$ is 
a $TK_5$ in  $G'$ with branch vertices $x_1,x_2,y_1,y_2,z_1$.  \qed

\end{document}